\documentclass[a4paper]{article}                    

\usepackage[margin=1.25in]{geometry}

\usepackage{graphicx}
\usepackage{amsmath}
\usepackage{amssymb}
\usepackage{amsfonts}
\usepackage{dsfont}
\usepackage{mathtools}
\usepackage{authblk}

\newcommand{\E}{\mathbb{E}}

\newcommand{\lam}{\lambda}
\newcommand{\m}{\mu}

\newcommand{\dt}{\delta t}
\newcommand{\intt}{\int\limits_{0}^{t}}
\newcommand{\etal}{\textit{et~al.}~}
\newcommand{\suml}{\sum\limits}

\newtheorem{theorem}{Theorem}[section]
\newtheorem{lemma}{Lemma}[section]
\newtheorem{remark}{Remark}[section]

\title{Job Allocation in Large-Scale Service Systems with Affinity Relations}

\author[1]{Ellen~Cardinaels\thanks{Corresponding author: Ellen~Cardinaels (e.cardinaels@tue.nl)}}
\author[1,2]{Sem~C.~Borst}
\author[1]{Johan~S.H.~van~Leeuwaarden}
\affil[1]{Eindhoven University of Technology, The Netherlands}
\affil[2]{Nokia Bell Labs, Murray Hill, USA}

\newcommand{\keywords}[1]{\textbf{Keywords} #1}
\newcommand{\MSC}[1]{\textbf{MSC2010} #1}

\begin{document}

\maketitle


\begin{center}
\maketitle
\end{center}

\begin{abstract} 
We consider load balancing in service systems with affinity relations between jobs and servers. Specifically, an arriving job can be allocated to a fast, primary server from a particular selection associated with this job or to a secondary server to be processed at a slower rate. Such job-server affinity relations can model network topologies based on geographical proximity, or data locality in cloud scenarios. We introduce load balancing schemes that allocate jobs to  primary servers if available, and otherwise to secondary servers. A novel coupling construction is developed to obtain stability conditions and performance bounds using a coupling technique. We also conduct a fluid limit analysis for symmetric model instances, which reveals a delicate interplay between the model parameters and load balancing performance.\\

\noindent \keywords{load balancing, stochastic coupling, fluid limit, job scheduling, network~topology}\\
\MSC{60K25, 68M20, 90B15, 90B22, 90B35}

\end{abstract}

\section{Introduction}
In this paper we analyze a load balancing scheme in a service system where particular servers are better equipped to process certain jobs because of affinity or compatibility relations. 
Load balancing algorithms play a crucial role in distributing jobs among multiple servers and have attracted strong renewed interest due to proliferation of large data centers and cloud computing. 
Well-known load balancing algorithms include for instance the Join-the-Shortest-Queue (JSQ), Join-the-Shortest-Queue-$d$ (JSQ($d$)) and Join-the-Idle-Queue (JIQ) policies.
These policies have been extensively analyzed in an overarching framework called the \textit{supermarket model}, consisting of a single dispatcher where jobs arrive that must be distributed among $N$ identical parallel servers.  
The JSQ policy assigns each arriving job to the server with the smallest queue length and has strong stochastic optimality properties among the class of policies without advance knowledge about the service requirements~\cite{ephremides1980simple,winston1977optimality}. The JSQ policy involves a significant communication burden however, which may be prohibitive in large systems.\\

This scalability issue has spurred an interest in the JSQ($d$) policy which assigns a job to the server with the smallest queue length among $d\ge 2$ randomly selected servers.
 Mitzenmacher~\cite{mitzenmacher2001power} and Vvedenskaya~\etal\cite{vvedenskaya1996queueing} analyzed the JSQ($d$) policy in an asymptotic regime where the total arrival rate and the number of servers grow large in proportion. Substantial performance gains were established compared to purely random assignment, even for $d=2$. Mukherjee~\etal \cite{mukherjee2018asymptotically} show that the waiting time in fact vanishes when $d$ tends to infinity as the number of servers grows large. A vanishing waiting time is also achieved by the JIQ policy which directs arriving jobs to an idle server or a randomly selected server if all servers are occupied \cite{lu2011join}. 
 The JIQ policy only has a constant communication overhead per job, but requires memory at the dispatcher. We refer to Van der Boor \etal \cite{van2017scalable} and Gamarnik \etal \cite{gamarnik2016delay} for further details.\\

A key feature of the supermarket framework is the exchangeability of the servers in the sense that any job can be handled equally well by any server, which is often not the case in practice. 
In the present paper we will focus on a scenario where jobs or servers are not intrinsically different, but where particular servers might be better equipped to process certain jobs because of affinity or compatibility relations.  
 Such affinity relations may for example arise due to geographical proximity in spatial settings, or data locality in content distribution or transaction processing applications.
The scenario will be modeled as follows: let $\mathcal{P}(\{1,\dots,N\})$ denote the power set of all servers. Then for a selection of servers, $S\in \mathcal{S} \subseteq \mathcal{P}(\{1,\dots,N\})$, jobs arrive %according to a Poisson process 
at rate $\lam_S \ge 0$. These jobs can be processed at rate $\m_1>0$ at any of the servers in $S$ or at rate $\m_2$ at any of the servers in $S^c$, with $\m_1>\m_2>0$. The arriving job is then labeled as a type~I or type~II job, depending on whether it can be served at rate $\mu_1$ or $\mu_2$, respectively. Our affinity-scheduling policy allocates the new job to a server in $S$ with the shortest queue length unless it might be beneficial to redirect the job to a server outside $S$. The precise allocation and scheduling strategies will be described in Section~\ref{subsec:modeldescr}. \\

When $\mathcal{S}$ contains all neighborhood sets of a graph $G_N$ on $N$ vertices, we refer to our model as the \textit{graph model}. 
The graph model extends the models constructed by Gast~\cite{gast2015power}, Turner~\cite{turner1998effect} and Mukherjee~\etal\cite{mukherjee2018asymptotically}. In these settings it is assumed that all nodes have equal arrival rates and jobs can only be forwarded to direct neighbors; it is not possible to redirect an arriving job to any other nodes.
The model constructed by Yekkehkhany~\etal \cite{yekkehkhany2018gb} does allow for jobs to be redirected to higher-degree neighbors to be served at lower rates. 
When $\mathcal{S}$ consists of all subsets of servers of fixed size $d$, we refer to our model as the \textit{combinatorial model}. In addition, the arrival rates will be equal among the server selections which strengthens the symmetric nature of the combinatorial model.  \\ 
 
 The lack of exchangeability among the servers makes the affinity-scheduling model complicated to analyze in general. The analytical techniques that are most commonly used in the context of the supermarket model, such as mean-field limits and even standard coupling arguments, fundamentally rely on this feature. 
 These techniques can only be applied for the combinatorial model. 
 For the general model, and in particular the graph model, the investigation of load balancing issues is challenging, and enters largely uncharted methodological territory.\\

We will establish a stochastic dominance result for the occupancy process of the general affinity-scheduling model, which will yield a sufficient stability constraint as an immediate by-product. 
Exploiting the coupling of this dominance result, we can derive two stronger dominance results for the graph model, which will in particular hold if the underlying graph structure is rather dense. To the best of our knowledge, these are the first results that explicitly capture the impact of network structure on load balancing performance.\\

 For the combinatorial model we will conduct a fluid-limit analysis. A trajectory of the fluid limit will converge to one of the possibly multiple fixed points, depending on the mutual relationships of the model parameters and the initial configuration of the system. When the fixed point is unique, we demonstrate that this provides a good approximation for the intractable stationary distribution in a finite server setting. When multiple fixed points occur, we observe the phenomenon of `tunneling' described by  \cite{gibbens1990bistability}. The stochastic process will switch between multiple modes corresponding to the locally stable fixed points of the fluid limit.\\
 
The remainder of this paper is organized as follows. A detailed model description will be provided in Section~\ref{subsec:modeldescr}. Next, the main stochastic dominance result is presented in Section~\ref{subsec:stochdom} together with the coupling that establishes this result. This section is completed with two stronger stochastic dominance properties for graph models. In Section~\ref{subsec:fluidlim} we present a fluid-limit analysis of the affinity-scheduling policy for combinatorial models. The proofs of all results are deferred to Section~\ref{sec:proofs}. Finally, Section~\ref{sec:conclusion} provides concluding remarks and some directions for further research.

%%%%%%%%%%%%%%%%%%%%%%%%%%%%%%%%%%%%%%%%%%%%%%%%%%%%%%%%%%%%%%%%%%%%%%%%%%%%%%%%%%%%%%%%%%%%%%%%%%%%%%%%%%%%%%%%%%%%%%%%%%%%%%%%%%%%%%%%%%%%%%%%%%%%%%%%%%
\section{Model description} \label{subsec:modeldescr}

We now describe the affinity-scheduling model with $N$ servers. 
For a selection $S\in \mathcal{P}(\{1,\dots,N\})$ jobs arrive as a Poisson process of rate $\lam_S \ge 0$. For these jobs, the servers in $S$ and $S^c$ are called the primary and secondary servers, respectively. An arriving job can be allocated as a type~I job to a primary server or as a  type~II job to a secondary server. Type~I jobs have independent and exponentially distributed service times with parameter $\m_1$. Type~II jobs have on average longer service times which are independent and exponentially distributed with parameter $\m_2 <\m_1$. It is important to note that the job type~is not predetermined on arrival, but established by the allocation strategy.
The main idea behind our allocation strategy is: `Allocate a job to a server in the primary selection unless it might be beneficial to allocate a job to a secondary server even though the service time might be longer'. The rationale for this is to reduce the waiting time of a job. More precisely, the \textit{allocation strategy} goes through the following three steps:
\begin{enumerate}
\item Is there at least one completely idle server in the primary selection~$S$? If so, allocate the arriving job as a  type~I job to one of these servers.
\item Is there at least one completely idle server in the secondary selection~$S^c$? If so, allocate the arriving job as a type~II job to one of these servers.
\item If there are no idle servers present, then allocate the job as a type~I job to the primary server with the smallest number of type~I jobs. Ties are broken according to the number of type~II jobs, in favor of a lower number.
\end{enumerate} 
When the second step is omitted, our policy resembles a JSQ($|S|$) policy with $|S|$ the cardinality of the primary selection. However, the cardinality of the server selection is allowed to differ among arriving jobs in our model and the server selection~$S$ itself is not sampled uniformly at random as is the case in a JSQ($|S|$) policy.
  Moreover, the second step can be related to a JIQ policy on the set of secondary servers. Our affinity-scheduling policy thus shares similarities with both policies.
We assume the time it takes for the dispatcher to make a decision for an arriving job, just as the possible time it takes a job to reach its selected server, to be negligible. Notice that due to this strategy, an arriving job will never be assigned as a type~II job to a server that already has a job in its queue. 
  Denote the configuration of a server, i.e. the number of type~I and type~II jobs in its queue, by $(i,j)$, $i,j\ge 0$. As an illustrative example of the allocation strategy, the primary and secondary servers have $\{(1,0),(1,1),(1,1),(4,0) \}$ and $\{(1,0),(1,1),(1,1),(3,1)\}$ as their configurations, respectively. Due to the allocation strategy, the third step will be applied and the primary server with configuration $(1,0)$ will receive an additional type~I job.
In general, type~I jobs are the preferred type of jobs, which also manifests itself in the \textit{scheduling strategy}. Each server operates under a preemptive priority scheduling discipline in favor of the type~I jobs. Moreover, type~I jobs are served in order of arrival.\\

Let $N_{n,j}(t)$ denote the number of type~$j$ jobs at server $n\in \{1,\dots,N\}$ at time $t$. The configuration of server $n$ is then given by $\left(N_{n,\text{I}}(t),N_{n,\text{II}}(t)\right) \in \mathbb{N}^2$ with state space $\mathbb{N}^{2N}$. The vector $\left(N_{n,\text{I}}(t),N_{n,\text{II}}(t)\right)_{n}$ evolves as an irreducible, time-homogeneous Markov process.
We also introduce different variables that are more server-centric and will be more convenient in proving stochastic dominance and analyzing the fluid limit. Define $Q^N_{ij}(t)$ as the number of servers with $i$ type~I jobs and $j$ type~II jobs at time $t$, with $i,j\ge0$. Then 
\begin{equation}
\overline{Q}^N_{ij}(t)\doteq \suml_{k\ge i} Q^N_{kj}(t)
\end{equation}
denotes the number of servers with at least $i$ type~I jobs and exactly $j$ type~II jobs. We note that $\sum_{j\ge0}\overline{Q}^N_{0j}(t)=N$ by definition. 
It is important to note that these variables will no longer lead to a Markov process representation in the general settings mentioned in the introduction. This immediately limits the number of available techniques to analyze the performance. \\

Let $\mathcal{S}$ denote the subset of $\mathcal{P}(\{1,\dots,N\})$ with strictly positive arrival rates. Besides the general setting where $\mathcal{S}$ can be any subset of $\mathcal{P}(\{1,\dots,N\})$ we will also investigate some more restricted cases. In the \textit{graph model} on graph topology~$G_N$, each node represents a server and the edges represent underlying relations between them. Then each set in $\mathcal{S}$ consists of a server and its neighbors determined by $G_N$. In total $\mathcal{S}$ contains $N$ sets and jobs arrive to each of these sets independently at a uniform rate $\lam$. 
 This setting mimics a situation where a job's physical arrival location plays a role in its service time at the various servers.\\ 

Let $\mathcal{S}$ consist of all possible server selections of size~$d$. The cardinality of $\mathcal{S}$ is $\binom{N}{d}$ and henceforth we refer to this model as the \textit{combinatorial model}. We assume a uniform arrival rate $\nu$ per selection. We let $\nu=\lam N / \binom{N}{d}$ per selection such that the total rate is given by $\lam N$. 
Observe that the combinatorial model captures the situation where a selection of $d$ servers is drawn uniformly at random as the primary selection for each job and arrivals occur at rate $\lam N$ in total.

\begin{remark}
There are also instances of the affinity-scheduling model that are not captured by either the graph model or the combinatorial model.
As an example, suppose a job arrives to a primary selection that consists of the servers $1,\dots,5$ or an arbitrary selection of size two of the remaining servers. Then $\mathcal{S}$ consists of $\{1,\dots,5\}$ and all pairs of  servers of $6,\dots,N$. For some $\nu>0$, the arrival rates per selection are given  by $\nu/2$ and $\nu/\left((N-5)(N-6)\right)$, respectively.
\end{remark}

%%%%%%%%%%%%%%%%%%%%%%%%%%%%%%%%%%%%%%%%%%%%%%%%%%%%%%%%%%%%%%%%%%%%%%%%%%%%%%%%%%%%%%%%%%%%%%%%%%%%%%%%%%%%%%%%%%%%%%%%%%%%%%%%%%%%%%%%%%%%%%%%%%%%%%%%%%
\section{Stochastic dominance and coupling} \label{subsec:stochdom}

In this section we establish several stochastic dominance results
for our affinity-scheduling strategy.
We will construct a stochastic coupling that allows a comparison
with various reference systems in terms of the ordered server states,
and refer to this coupling as the \textit{affinity coupling}. In contrast to the original system with affinity relations,
the various reference systems all involve $N$~exchangeable servers,
and are therefore far more amenable to (asymptotic) analysis,
yielding tractable performance bounds.
We will not explicitly consider the type~II jobs since our allocation
strategy will never add such a job at any server where there are already type~II jobs present and we assume that the initial configuration of the system will only have a finite number of type~II jobs.
Thus we focus on the type~I jobs that enjoy higher affinity at each
of the servers, hence the name of the coupling. \\

While each of the $N$~servers in the reference system processes jobs
in a FCFS manner at rate~$\mu_1$, the various specific incarnations
differ in the value of the normalized arrival rate per server~$\lambda_0$
and the policy for assigning jobs.
The choice of the specific reference system is aligned with the properties of the
original system in terms of the server selections~${\mathcal S}$ and the
associated arrival rates $\lambda_S$, $S \in {\mathcal S}$.
Loosely speaking, we obtain increasingly strong dominance results under
increasingly restrictive symmetry and structural conditions on the
server selections~${\mathcal S}$ and the associated arrival rates.
The three specific variants for the reference system that we consider
operate under either
(i) a purely random assignment (RA) policy,
(ii) a MJSQ($k$) policy (as specified later), or
(iii) a JSQ($k$) policy (as described in the introduction).
While the RA system provides exact upper bounds in terms of independent
M/M/1 queues, the MJSQ($k$) and JSQ($k$) systems yield asymptotic
upper bounds based on fluid limits.

The dominance results revolve around stochastic majorization properties
in terms of the ordered server states.
Specifically, define
\begin{equation}
\overline{Q}_{m+}^{\mathrm{aff}}(t) \doteq \sum_{i = m}^{\infty} \overline{Q}_i^{\mathrm{aff}}(t) \text{~and~} \overline{Q}_{m+}^{\mathrm{ref}}(t) \doteq \sum_{i = m}^{\infty} \overline{Q}_i^{\mathrm{ref}}(t),
\end{equation}
with $\overline{Q}_i^{\mathrm{aff}}(t)$ and $\overline{Q}_i^{\mathrm{ref}}(t)$ denoting the number of servers with at least $i$ type~I jobs in their queue at time~$t$ in the original system and the reference system, respectively.
We will establish results of the form
\[
\{(\overline{Q}_{m+}^{\mathrm{aff}}(t))_{m \geq 1}\}_{t \geq 0} \leq_{{\mathrm{st}}}
\{(\overline{Q}_{m+}^{\mathrm{ref}}(t))_{m \geq 1}\}_{t \geq 0}.
\]
This majorization result indicates that the number of type~I jobs
residing in the $m$-th or higher queue position in the original system
is stochastically bounded from above by the number of jobs residing
in the $m$-th or higher queue position in the reference system.
In particular, taking $m = 1$, this implies that the total number
of type~I jobs in the original system is stochastically bounded from above
by the total number of jobs in the reference system.
As noted earlier, we know the exact distribution of the latter quantity
in the RA system and have an asymptotic result for the MJSQ($k$)
and JSQ($k$) systems.

\begin{figure}[h]
\centering
\includegraphics[width=13cm]{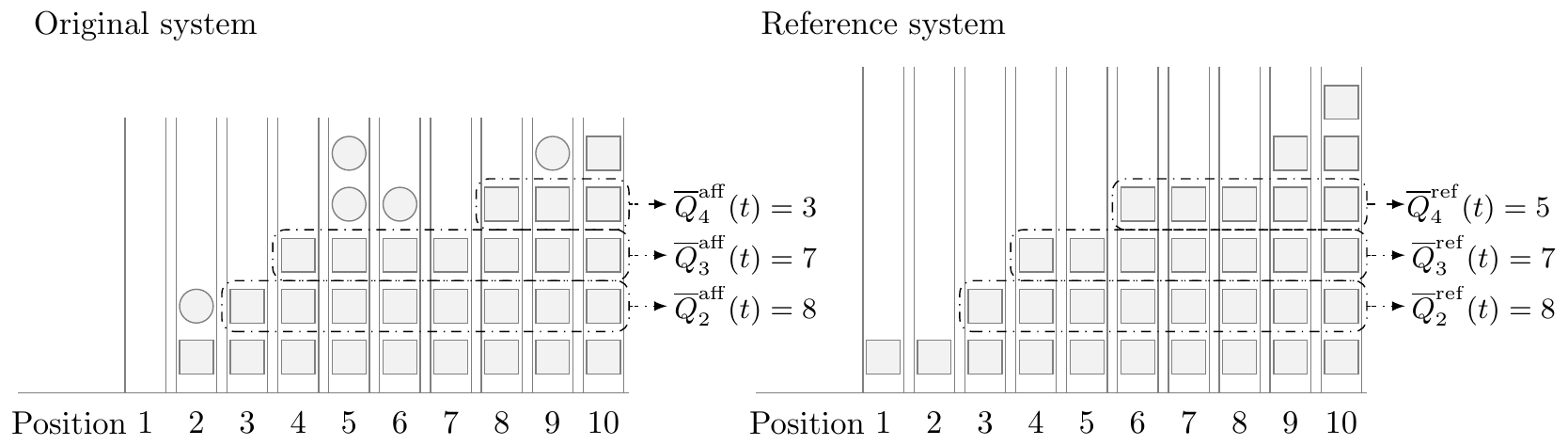}
\caption{A schematic representation of the two ordered systems with $N=10$ servers at time $t$. The original system, operating under the affinity-scheduling policy, consists of type~I (squares) and type~II (circles) jobs. The $N$~exchangeable servers in the reference system make no distinction between the jobs. } \label{fig:sample_path}
\end{figure}

In order to prove the stochastic majorization properties, we introduce
the \textit{affinity coupling} to construct sample paths for the
original and reference systems on a joint probability space for which
the stated inequalities hold in a deterministic way
\cite{liu1995sample,sparaggis1994sample,stoyan1983comparison}. 
For all three specific reference systems, the common proof concept is
to ensure that under the coupling two key properties always hold
with respect to the ordered server states as illustrated in Figure~\ref{fig:sample_path}:
 (a) addition of a type~I job at an arrival epoch in the
original system must be accompanied by insertion of a job
at a higher-ordered server in the reference system; (b) removal of a job at a service completion epoch in the reference
system must force disposal of a type~I job at the same ordered server
in the original system (unless there is no type~I job at that server). We can prove the following general lemma.

\begin{lemma}\label{lem:generalcoupling} \emph{(Affinity coupling).} If a stochastic coupling between the original system and the reference system can be constructed such that $(a)$ and $(b)$ are satisfied, then $(\overline{Q}_{i}^{\mathrm{aff}}(t))_{i\ge 1}$ is majorized by $(\overline{Q}_{i}^{\mathrm{ref}}(t))_{i\ge 1}$ for $t\ge 0$. In the sense that
\begin{equation} \label{eq:lem_maj}
\suml_{i = m}^{ \infty}\overline{Q}_{i}^{\mathrm{aff}}(t) \le \suml_{i = m}^{ \infty}\overline{Q}_{i}^{\mathrm{ref}}(t)
\end{equation} 
for all $m\ge 1$, provided that the initial configurations of both systems satisfy this inequality.
\end{lemma}
The proof of Lemma~\ref{lem:generalcoupling} can be found in Subsection~\ref{subsec:proof_stoch_dom}.
In the remainder of this section, we will precisely describe the affinity coupling for each of the reference systems under consideration and verify that the properties (a) and (b) are satisfied. The coupling at service completion epochs to ensure property~(b) as
further detailed below is fairly standard and common across all three
reference systems.\\

\noindent \textbf{Coupling at arrival epochs.}
In contrast to the service completions, the coupling at arrival epochs to guarantee property~(a) is
novel and highly specific to the reference system under consideration.
Due to the lack of exchangeability among servers, the coupling
at arrival epochs involves a further subtle complication that does not
arise in constructing sample path comparisons in the context of the
ordinary supermarket model.
Even though we compare the evolution of the two systems in terms of the
$\overline{Q}_i$ variables as usual, these generally do not provide
a Markovian state description for the original system as noted earlier in Section~\ref{subsec:modeldescr}.
In particular, the transitions at arrival epochs intricately depend
on the server selections~${\mathcal S}$, and cannot be suitably represented
in terms of the $\overline{Q}_i$ variables. \\

\noindent \textbf{Coupling at service completion epochs.}
The coupling generates potential service completions at rate $\m_1 N$, but the aggregate service rate in either the original or the reference system might be lower than $\m_1N$ because of servers being idle or only working at rate $\m_2$ on type~II jobs.
 Let $W_{\mathrm{aff}}$ and $W_{\mathrm{ref}}$ be the sets of ordered positions of servers in the original and reference system, respectively, that are working on (type~I) jobs just before some time $t$ at which a potential service completion occurs. Define $W$ as the intersection $W_{\mathrm{aff}} \cap W_{\mathrm{ref}}$ which equals $W_{\mathrm{aff}}$ or $W_{\mathrm{ref}}$ due to the ordering and the preemptive strategy of the affinity-scheduling policy.
 A random variable $X_t$, drawn from a uniform distribution on $[0,1]$, decides which of the following actions is selected. 
\begin{itemize}
\item[(i)] $0\le X_t \le \frac{|W|}{N}$: Sample uniformly at random a position $n$ from $W$; a departure will take place at time $t$ in both the original and the reference system at the server located at position $n$.
\item[(ii)] $\frac{|W|}{N} < X_t \le \frac{|W_P|}{N}$ where $P$ is `$\mathrm{aff}$' or `$\mathrm{ref}$': Sample uniformly at random a server position from $W_P\setminus W$; one job will be removed from the corresponding server in system $P$ at time $t$.
\item[(iii)] $X_t  > \frac{\max\{|W_{\mathrm{aff}}|,|W_{\mathrm{ref}}|\}}{N}$: No real departure will occur among the type~I jobs in the original system or the jobs in the reference system.
\end{itemize}
We note that the total departure rate of type~I jobs from the original system is indeed given by  $\m_1 |W_{\mathrm{aff}}|$, likewise for the reference system with a total departure rate of $\m_1 |W_{\mathrm{ref}}|$. The idea to work with intersections of the active server sets comes from~\cite[Section 4]{mukherjee2016universality}.

\subsection{Affinity coupling with the general model}

We now consider a general structure for the server selections~$\mathcal{S}$ and the corresponding arrival rates $\{\lam_S \mid S \in  \mathcal{S}\}$ per server selection. The reference system will operate under the RA policy with arrival rate $\lam_0$ per server. Thus $\lam_0<\m_1$ is a sufficient stability condition for the reference system. So the purpose of this subsection is twofold: the affinity coupling is illustrated in a general setting of our affinity-scheduling policy in order to obtain a stochastic dominance result and a stability condition is obtained as an immediate by-product. \\

The choice of $\lambda_0$ is determined by the arrival rates per server selection in the original system, namely
\begin{equation}\label{eq:lam0}
\lam_0\doteq \underset{p_{Sn}}{\min}\left\{ \underset{n}{\max} \left\{ \lam_n^* = \sum_{S\in\mathcal{S}: n\in S}\lam_{S}p_{Sn} \mid \sum_{n\in S} p_{Sn} = 1 \text{~with~} p_{Sn} \ge 0, \forall n\in S \right\}\right\}.
\end{equation}
The variable $p_{Sn}$ may be interpreted as the fraction of jobs with server selection~$S$ that are assigned to server~$n\in S$. With this interpretation in mind, it is easily seen that at least one server must handle an arrival rate of $\lam_0$ or larger in case jobs are only allowed to be executed as type~I jobs. Thus $\lam_0<\m_1$ is clearly a necessary stability condition for any strategy in this case. The condition is sufficient as well, for instance for a simple static strategy that assigns a job with server selection~$S$ to server $n$ with probability $p_{Sn}$. However, the implementation of this strategy would require full knowledge  of the arrival rates $\lam_S$. We will establish  that the condition is also sufficient for the stability of our affinity-scheduling strategy, which does not rely on any knowledge of the arrival rates $\lam_S$ at all.\\

We now specify the affinity coupling for the reference system with the RA policy.\\

\noindent \textbf{Coupling at arrival epochs.} The coupling generates potential arrival events at rate $\lam_0 N$. If a potential arrival occurs at time $t$, a position $n^*$ from the set $\{1,\dots,N\}$ is selected uniformly at random. For brevity we simply refer to the server at position $n^*$ as server $n^*$. An addition  of a new job in the reference system will take place at this server~$n^*$. Since this position was randomly selected, the coupling strategy will give rise to an addition according to the RA policy in a system with arrival rate $\lambda_0$ per server.

 In order to determine whether an arrival event of a type~I job takes place in the original system and at which server this will happen, we follow the below-described strategy.
 Two random variables, $Y_{t,1}$ and $Y_{t,2}$, are sampled from a uniform distribution on $[0,1]$ to take into account that the total arrival rate in the original system might be smaller than $\lam_0 N$ and to select a server selection~$S$ for an arriving job. To make the decisions, we rely on the variables $(p_{Sn}^*)_{S,n}$ that attain the minimum in~(\ref{eq:lam0}).
First, $Y_{t,1}$ establishes if an arrival occurs to a primary selection containing server $n^*$, which happens with probability $\lam_{n^*}^*/\lam_0$. If an arrival will take place, then a server selection~$S$ containing $n^*$ is selected as the primary selection with probability $\lam_{S}p_{Sn^*}^*/\lam^*_{n^*}$ for which $Y_{t,2}$ is used. All remaining servers form the secondary selection. Note that the total arrival rate to a server selection~$S$ 
\begin{equation}
\lam_0 N \sum_{n\in S} \frac{1}{N}\frac{\lam^*_n}{\lam_0}\frac{\lam_{S}p_{Sn}^*}{\lam^*_n}
\end{equation}
will indeed be equal to $\lam_S$ in the original system as $\sum_{n\in S} p_{Sn}^*=1$ by definition such that this method to handle arriving jobs will coincide with the arrival process of general model described in Section~\ref{subsec:modeldescr}. Once these selections are set for an arriving job, we apply the allocation policy as defined in Section \ref{subsec:modeldescr}. Due to the general structure of $\mathcal{S}$ it is not possible to determine the exact server at which a job is allocated in terms of the variables $(\overline{Q}^N_{ij})_{i,j}$. However, if the new job is allocated as a type~I job in the original system to one of the servers in $S$, it is known that the position of this server will be at most~$n^*$. Since the newly arrived job in the reference system is assigned to server~$n^*$, property~(a) of the coupling is maintained.\\

With the notation as introduced above, we can prove the following theorem.

\begin{theorem}\label{th:general_dom}\emph{(General affinity-scheduling model).}
Let $\lambda_0$, as defined in~$\mathrm{(\ref{eq:lam0})}$, be the arrival rate per server in the reference system operating under the \emph{RA} policy. Then, for suitable initial conditions,
\begin{equation}
\{(\overline{Q}_{m+}^{\mathrm{aff}}(t))_{m \geq 1}\}_{t \geq 0} \leq_{{\mathrm{st}}}
\{(\overline{Q}_{m+}^{\mathrm{RA}}(t))_{m \geq 1}\}_{t \geq 0}
\end{equation}
holds for the general affinity-scheduling model with $N$~servers.
\end{theorem}

The above-described coupling between the original and the reference system operating under the RA policy satisfies the general framework of the affinity coupling, i.e. properties (a) and (b) hold. Then the result stated in Lemma~\ref{lem:generalcoupling} is applicable so that Theorem~\ref{th:general_dom} follows from the majorization result established there. 
Theorem~\ref{th:general_dom} provides a stochastic upper bound for the total number of type~I jobs in the original system in terms of the number of jobs in a reference system with the RA policy by taking $m=1$. Although this upper bound is sufficient to guarantee stochastic stability for $\lam_0<\m_1$, we will develop stronger majorization results for particular settings of the graph model in the next two subsections. The method to prove the two stronger results is captured in the general framework stated at the beginning of this section, but requires a different and more advanced coupling method between arrival events. 

\subsection{Graph model}

We will further investigate our model on a graph topology $G_N$ as described in Section~\ref{subsec:modeldescr}. It is challenging to get a grip on the performance of an allocation policy that is applied in a network structure, and establishing stochastic dominance relations can give an initial insight into the theoretical behavior of load balancing algorithms in structured environments. It was mentioned in  Section~\ref{subsec:modeldescr} that the arrival rate over all server selections established by the graph structure $G_N$ is given by $\lam$, and thus Theorem~\ref{th:general_dom} is still valid if we set $\lambda_0 = \lambda$.
However, we will make two different assumptions on the structure of the graph topology and for each of them a much stronger dominance result than Theorem~\ref{th:general_dom} is obtained. The first scenario assumes that the minimum degree of $G_N$ is sufficiently high and the second scenario entails regular graph topologies.

\subsubsection{Minimum degree}

The reference system with $N$~exchangeable servers operates under an allocation policy related to JSQ, namely MJSQ($k$) \cite{mukherjee2018asymptotically}. 
In this setting new jobs arrive at a total rate of $\lam N$ and are processed at a server according to a FCFS policy at rate $\m_1>\lam$. An arriving job is allocated to the server with the $(k+1)$-th smallest queue length. A clear analogy can be seen if the system is initially completely empty; then $k$ servers will constantly remain idle. The system operates as if only $N{-}k$ servers are present and applies a JSQ policy restricted to these servers. If $N$ is sufficiently large compared to $k$, i.e. if
\begin{equation}\label{eq:restr_MJSQ}
\lam N < \m_1 (N-k),
\end{equation}
then the MJSQ($k$) policy is stochastically stable. It is intuitively clear that this policy can achieve much better performance than the RA policy.\\

Suppose that the minimum degree of the graph $G_N$ is at least $N{-}k{-}1$, without any other structural assumptions. Let $N$ and $k$ satisfy the relation in~(\ref{eq:restr_MJSQ}), then we can describe a coupling between our graph model with underlying topology~$G_N$ and the reference system with the MJSQ($k$) policy. The coupling between both systems will fit the general framework of the affinity coupling but the coupling method for the arriving jobs will differ from the general setting in the previous subsection.\\

\noindent \textbf{Coupling at arrival epochs.} For each of the neighborhood sets in $\mathcal{S}$ there is a uniform arrival rate $\lam$ such that the total arrival rate in the original system is also given by $\lam N$. Assuming that an event in the coupled sample path is an arrival, it is always directed to the server at position $k+1$ under the MJSQ($k$) policy. For the graph model, the primary selection~$S$ consists of a randomly selected server and its neighbors under the topology~$G_N$ and the secondary selection~$S^c$ contains all other servers. We do not know the exact ordered positions of the servers in the primary selection that is of size at least $N-k$ in terms of the $\overline{Q}_i$ variables. The worst-case scenario that could arise is a primary selection of size exactly $N{-}k$ where the servers are the $N{-}k$ highest ordered servers. Then a type~I job is allocated to the server at position $k+1$. All other scenarios where an arriving job is labeled as a type~I job in the original system will lead to an allocation that is at most at the $(k+1)$-th position. Hence property (a) of the affinity coupling is satisfied.\\

Together with the coupling between service completions as explained in the first part of this section, the coupling between our affinity-scheduling policy on a graph structure with minimum degree $N{-}k{-}1$ and the reference system with the MJSQ($k$) policy satisfies the general framework of the affinity coupling stated in Lemma~\ref{lem:generalcoupling}. Thus Theorem~\ref{th:maj_min_deg} follows from the majorization result in Lemma~\ref{lem:generalcoupling}.

\begin{theorem}\label{th:maj_min_deg}\emph{(Graph model with minimum degree $N{-}k{-}1$).}
Consider the graph model with an underlying graph topology with minimum degree $N{-}k{-}1$ and a reference system that operates under the \emph{MJSQ($k$)} policy. Then, for suitable initial conditions, 
\begin{equation}
\{(\overline{Q}_{m+}^{\mathrm{aff}}(t))_{m \geq 1}\}_{t \geq 0} \leq_{{\mathrm{st}}}
\{(\overline{Q}_{m+}^{{\small \text{\emph{MJSQ($k$)}}}}(t))_{m \geq 1}\}_{t \geq 0}.
\end{equation}
\end{theorem}

Once the reference system is stochastically stable, if condition~(\ref{eq:restr_MJSQ}) is fulfilled, we can give a meaningful upper bound on the total number of type~I jobs in the graph model in terms of the total number of jobs under the MJSQ($k$) policy. This upper bound will be stronger compared to the result in Theorem~\ref{th:general_dom} since the MJSQ($k$) policy outperforms the RA policy with arrival rate $\lambda$ per server and service rate $\m_1$. 

\begin{remark}
Theorem~\ref{th:maj_min_deg} can be generalized for scenarios where each server selection~$S$ has a size of at least $N{-}k$ and a non-uniform arrival rate $\lam_S$.
\end{remark}

\subsubsection{Regular graph}

As mentioned in the introduction, JSQ($k$) gives already substantial performance improvements for small values of $k$ compared to the RA policy. With this in mind, we show that the number of type~I jobs under our affinity-scheduling policy on a $d$-regular graph is stochastically dominated by the total number of jobs under a JSQ($k$) policy, when $d$ and $k$ satisfy the following relation:
\begin{equation}\label{eq:graph_mayo_reg}
\sum_{i=1}^{N-d-1} \binom{N-i}{k-1} < \frac{d+1}{N} \binom{N}{k}.
\end{equation}

The proof requires a coupling between the arrival events in both systems such that feature (a) of the affinity coupling is maintained. 
 We will introduce a novel approach to represent or visualize all possible server selections in $\mathcal{S}$ that an arriving job can choose from. \\ 

\begin{figure}[h]
\centering
\includegraphics[width=9cm]{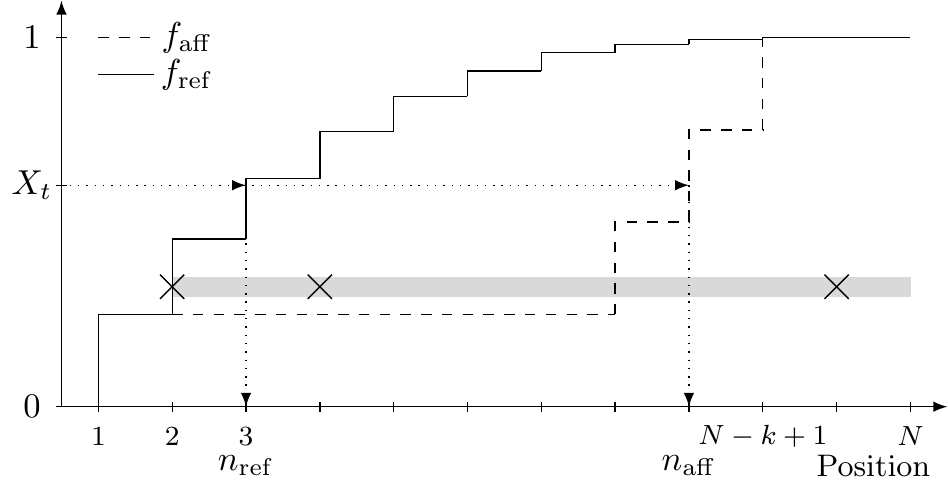}
\caption{A visualization of the step functions of the server selections in both systems. There are $N=12$ servers, $d=2$ and $k=3$. The block interpretation of server selection~$S=\{2,4,11\}$ of the reference system is represented in gray.} \label{fig:stepfunction}
\end{figure}

An arriving job in the reference system with JSQ($k$) is allocated to the lowest positioned server among $k$ randomly selected servers. In total there are $\binom{N}{k}$ server selections and each server belongs to $\binom{N-1}{k-1}$ different server selections. Thus, the lowest positioned server of the system belongs to $\binom{N-1}{k-1}$ different server selections, the second lowest server is part of precisely $\binom{N-2}{k-1}$ different server selections without the lowest ordered server. One can continue this reasoning up to the $(N-k+1)$-th lowest ordered server: this server belongs to only one more server selection that is not yet observed at any of the lower ordered servers. All higher ordered servers cannot be part of an unobserved server set. We will construct a step function from the positions $\{1,\dots,N\}$ to the interval $[0,1]$ in order to represent the server selections. 
Assume that the servers are ordered from 1 to $N$ and the lowest ordered position of $k$~selected servers is denoted by $n$. We represent this selection as a block from position $n$ to $N$ with height $1/\binom{N}{k}$. This procedure can be repeated for each of the $\binom{N}{k}$ possible selections. Stacking all these blocks according to their length will give rise to the following step function: 
\begin{equation}
f_{\mathrm{ref}} \colon \{1,\dots,N\} \to [0,1] \colon x \mapsto \left\{
\begin{array}{lcr}
\frac{1}{\binom{N}{k}}\sum_{i=1}^{x} \binom{N-i}{k-1}, & ~&1\le x\le N-k+1,\\
1, & ~& N-k+1 < x \le N.
\end{array}
\right.
\end{equation}
An example of this visualization can be found in Figure \ref{fig:stepfunction}.
%%%%%%%%%%%%%%%%%%%%%%%%%%%%%%%%%%%%%%%%%%%%%%%%%%%%%%%%%%%%%%%%%%%%
%%%%%%%%%%%%%%%%%%%%%%%%%%%%%%%%%%%%%%%%%%%%%%%%%%%%%%%%%%%%%%%%%%%%%%%

We aim to construct a similar step function based on the possible primary server selections for the original system when the underlying graph topology is a $d$-regular graph. Jobs arrive at a total rate $\lam N$ and an arriving job selects uniformly at random a server selection~$S$ from $\mathcal{S}$. By construction, $\mathcal{S}$ contains $N$ different primary server selections, each of size $d+1$. Then, the lowest ordered server in the system belongs to $d+1$ different server selections. However, it is not possible to count the number of additional server selections containing the second lowest ordered server without knowing the position of each of the servers, since we operate on a fixed graph structure. We construct a step function based on the worst-case scenario where the lowest positioned server of each of the server selections is still at the highest possible position. The first jump occurs at the lowest positioned server, while all remaining jumps will occur at the highest possible positioned servers. This induces stronger correlations between the servers that have more type~I jobs. Notice that the worst-case ordering of the servers might not be a valid $d$-regular structure, so that this approach might be too conservative. The step function of this worst-case scenario is given by

\begin{equation}
\begin{array}{rl}
&f_{\mathrm{aff}}  \colon  \{1,\dots,N\} \to [0,1]\\
& x \mapsto \left\{
\begin{array}{lcr}
\frac{d+1}{N}, &~& 1\le x\le N-d - \lceil\frac{N}{d+1}\rceil+2\\
\frac{d+1}{N} +\frac{1}{N}\left(N-\lfloor\frac{N}{d+1}\rfloor (d+1)\right), &~& \frac{N}{d+1}\in \mathbb{N} \text{~and~} x = N-d - \lfloor\frac{N}{d+1}\rfloor+2 \\
1-\frac{d+1}{N}(N-d-x), &~& N-d - \lfloor\frac{N}{d+1}\rfloor+2 \le x < N-d\\
1, &~& N-d\le x \le N.
\end{array}
\right.
\end{array}
\end{equation}
An example of this step function can be found in Figure \ref{fig:stepfunction}.\\

\noindent \textbf{Coupling at arrival epochs.} 
Let the total arrival rate be $\lam N$ in both systems. 
For an arriving job at time $t$ we determine the servers of interest using the inverse transform sampling method \cite[Chapter~2]{devroye1986nonuniform}.
First, we note that the functions $f_{\mathrm{aff}}$ and $f_{\mathrm{ref}}$ are cumulative distribution functions by construction.
 Second, the only server of interest of the server selection $S$ in the original system or the server selection in the reference system is the lowest positioned server. So we sample a random variable $X_t$ from a uniform distribution on $[0,1]$ and determine the two servers positions, $n_{\mathrm{aff}}$ and  $n_{\mathrm{ref}}$, of interest of both systems. In the original system a server can be allocated as a type~I job to the selected server or to any other server as a type~II job. This procedure is visualized in Figure~\ref{fig:stepfunction}. 

So in order to guarantee feature (a) of the affinity coupling, it needs to be ensured that $n_{\mathrm{aff}}\le n_{\mathrm{ref}}$.
 Feature (a) is guaranteed if the step function of the original system is above the step function of the reference system, i.e. $f_{\mathrm{aff}}(n)\ge f_{\mathrm{ref}}(n)$ for all positions $n$. %  due to the interpretation of the blocks as server sets. 
 First we observe that the $d$-regular graph must be rather dense in order to obtain a stronger upper bound than provided by the RA policy. This can be seen if we investigate the step function at position 
 \begin{equation}
 x=N-d-\left(\left\lceil\frac{N}{d+1}\right\rceil-1\right).
 \end{equation}
  Once the degree $d$ is at least $N/2$, it is straightforward to show that $f_{\mathrm{aff}}(x)\ge f_{\mathrm{ref}}(x)$. This immediately implies that the step function $f_{\mathrm{aff}}$ only makes two jumps, at positions 1 and $N-d$ of sizes $(d+1)/N$ and $(N-d-1)/N$, respectively. Since the step function $f_{\mathrm{ref}}$ is concave in its discrete points, we only need to ensure that $f_{\mathrm{aff}}(N-d-1) > f_{\mathrm{ref}}(N-d-1)$ holds so that the step function of the original system is above the step function of reference system. This results in condition~(\ref{eq:graph_mayo_reg}) on the values of $d$ and $k$. Due to the coupling construction, we can prove the following dominance result.

\begin{theorem}\label{th:maj_dreg}\emph{(Graph model with $d$-regular graph).}
Consider the graph model with an underlying $d$-regular graph topology and a reference system operating under a \emph{JSQ($k$)} policy. If the model parameters $d$ and $k$ satisfy condition~$(\ref{eq:graph_mayo_reg})$, then, for suitable initial conditions,  
\begin{equation}
\{(\overline{Q}_{m+}^{\mathrm{aff}}(t))_{m \geq 1}\}_{t \geq 0} \leq_{{\mathrm{st}}}
\{(\overline{Q}_{m+}^{{\small \text{\emph{JSQ($k$)}}}}(t))_{m \geq 1}\}_{t \geq 0}.
\end{equation}
\end{theorem}

Due to the coupling construction using the block interpretation of the server selections, the above-described coupling fits the general framework of the affinity coupling as stated in Lemma~\ref{lem:generalcoupling}. Therefore, the result of Theorem~\ref{th:maj_dreg} follows from the result in Lemma~\ref{lem:generalcoupling}. If $\lam < \mu_1$, the reference system is stochastically stable and provides a meaningful upper bound on the performance of the graph model on a $d$-regular topology.

We list in Table~\ref{tab:grap_majo_reg} the minimum value of $d$ as a function of $k$ that guarantees the required dominance of the step functions for a system with $N=50$ servers.
 We observe that the graph structure is rather dense in order to stochastically dominate our process with a JSQ($k$) policy even for small values of $k$. One can argue that the primary selection of our affinity scheduling strategy must be much larger compared to the server selection under JSQ($k$) in order to guarantee better performance. But we should keep in mind that the underlying graph structure is fixed and all possible server selections are predetermined while the JSQ($k$) strategy can be seen as a strategy on a complete graph where an arbitrary set of size $k$ of the servers can be selected. This will affect the performance compared to a system with $N$ exchangeable servers, which is intuitively clear. 
\begin{table}[h]
\centering
\begin{tabular}{c|ccccccc}
$k$ & 2  & 3  & 4  & 5  & 10 & 15 & 25 \\ \hline
$d$ & 31 & 34 & 36 & 38 & 42 & 44 & 46
\end{tabular}
\caption{ The smallest possible value of $d$ that satisfies condition (\ref{eq:graph_mayo_reg}) is listed for a system with $N = 50$ servers and a given value of $k$.}
\label{tab:grap_majo_reg}
\end{table}
Moreover, it is important to note that the obtained value of $d$ might be too conservative. Our coupling method using the step functions requires a degree that is at least equal to $N/2$ in order to upper bound by the strategy JSQ(1), i.e. the  RA policy. On the other hand, we showed in the general result that the number of type~I jobs under any structural interpretation $\mathcal{S}$ is stochastically dominated by the number of jobs under a random assignment strategy.

\begin{remark}\label{rem:combmodel}
\emph{(Combinatorial model).}
Applying our affinity-scheduling strategy to the combinatorial model with $N$~servers shows a lot of similarities with the JSQ($d$) policy in a setting of $N$~exchangeable servers. Namely, an arriving job is allocated to the server with the shortest queue length among $d$ arbitrarily selected servers and sometimes the job can be directed to an idle server outside this selection. The coupling can be adjusted such that the number of type~I jobs at the $n$-th ordered position under our affinity-scheduling policy is less than or equal to the number of jobs at the $n$-th ordered position under a JSQ($d$) policy. The relation between the two policies will become more apparent in Section~\ref{subsec:fluidlim} when fluid-limit results are investigated.

Moreover, it can be shown that the combinatorial model is stochastically stable under a preemptive and a non-preemptive scheduling strategy using a Foster-Lyapunov argument. This result is shown under the assumption that the number of type~II jobs at each server never exceeds one. When the initial condition already satisfies this feature, the affinity-scheduling policy will never add a second type~II job to a server. The fact that stability is preserved under a preemptive strategy in favor of the type~I jobs is no surprise due to the structure of the server selections~$\mathcal{S}$ and the resemblance of the first step in the allocation strategy with the JSQ($d$) policy. Under a non-preemptive strategy it is no longer intuitively clear, as any finite value of $\m_2$ is allowed and one could imagine a situation where all servers are processing a type~II job and type~I jobs start to accumulate behind these type~II jobs.
\end{remark}

%%%%%%%%%%%%%%%%%%%%%%%%%%%%%%%%%%%%%%%%%%%%%%%%%%%%%%%%%%%%%%%%%%%%%%%%%%%%%%%%%%%%%%%%%%%%%%%%%%%%%%%%%%%%%%%%%%%%%%%%%%%%%%%%%%%%%%%%%%%%%%%%%%%%%%%%%%
\section{Fluid limit and fixed point analysis} \label{subsec:fluidlim}

As mentioned in the introduction, the affinity model in general lacks the exchangeability among the servers that underpins the use of mean field limits as the main analytical techniques in the supermarket model. 
Due to its inherent symmetry, the combinatorial model with uniform arrival rates for each of the server selections in $\mathcal{S}$ as described in Section~\ref{subsec:modeldescr} is one of the exceptions.
The variables $(Q_{ij}^N(t))_{i,j}$ will give rise to a Markov process representation in this case.
The primary and secondary server selections for an arriving job are of sizes $d_1$ and $d_2 = N-d_1$, respectively; we refer to both selections as the $d_1$-selection and $d_2$-selection. 
In order to gain insight in the system performance, we introduce the fluid scaled variables, i.e.
\begin{equation*}
\left(\frac{\bar{Q}^N_{i,j}(t)}{N}\right)_{i,j},
\end{equation*}
and analyze a sequence of systems where the number of servers~$N$ tends to infinity. The (weak) limit that arises is referred to as the fluid limit and is denoted by $(\overline{q}_{ij}(t))_{i,j}$.
When it is helpful to stress the proportion of servers with exactly $i$ type~I jobs, instead of at least $i$ type~I job, we consider the variables $(q_{ij}(t))_{i,j}$.
 In this section we consider initial configurations of the system that give rise to a process with at most one type~II job at each server through the complete process.
 as mentioned in Remark~\ref{rem:combmodel}. Furthermore, we assume that $\lam<\m_1$ to guarantee stochastic stability.
% {\color{blue}A characterization of the (deterministic) fluid limit will be provided in Subsection~\ref{subsec:fluidlimconstruction}. Next, the existence of its fixed points is shown in Subsection~\ref{subsec:fixedpoint}. Finally, the results are interpreted in Subsection~\ref{subsec:analysis}. }
Throughout this section we will consider a system with $\lam=0.8$, $\m_1=1$, and $\m_2=0.5$ in the numerical and simulation experiments, unless specified otherwise.

\subsection{Fluid limit}\label{subsec:fluidlimconstruction}
We now provide a characterization of the (deterministic) fluid limit in terms of a set of discontinuous differential equations. The $t$ reference in the notation will be omitted, if the context allows this.

We introduce a reduced arrival rate $\tilde{\lam}$. A job will always be directed to an idle server if available, either as a type~I job or a type~II job and idle servers are generated at rate $\m_1q_{10}-\m_2q_{01}$. This implies that if $\lam$ is sufficiently high, i.e. $\lam  > \m_1q_{10}+\m_2q_{01}$, only a fraction of the arriving jobs will start to queue in front of a server as type~I jobs on  fluid level. This fraction is given by $\tilde{\lam}/\lam$ with
\begin{equation}\label{eq:lam_tilde}
\tilde{\lam}=\left( \lam - \m_1q_{10}-\m_2q_{01}\right)^+.
\end{equation} Then,
\begin{equation} \label{eq:fluid_limit}
\begin{dcases}
\tfrac{d}{dt} \overline{q}_{00} = \m_2(\overline{q}_{01} - \overline{q}_{11}) - \lam(1-q_{00})^{d_1} + \tilde{\lam}  \mathds{1}\{q_{00} = 0\} \\
\tfrac{d}{dt} \overline{q}_{01} = \m_2(\overline{q}_{11}-\overline{q}_{01}) + \mathds{1}\{q_{00} > 0\} \left[ \lam(1-q_{00})^{d_1}\right] + \mathds{1}\{q_{00} = 0\} \left[ \lam - \tilde{\lam}\right]\\
\tfrac{d}{dt} \overline{q}_{10} =  \m_1(\overline{q}_{20} - \overline{q}_{10}) + \mathds{1}\{q_{00} > 0\}\left[\lam\left(1-(1-q_{00})^{d_1}\right)\right]\\
\tfrac{d}{dt} \overline{q}_{11} = \m_1(\overline{q}_{21} - \overline{q}_{11}) + \tilde{\lam } \mathds{1}\{q_{00} = 0\} \left[ \left(\overline{q}_{10}+\overline{q}_{01}\right)^{d_1} -\left(\overline{q}_{10}+\overline{q}_{11}\right)^{d_1}\right]\\
\text{for~}i\ge 2,\\
\tfrac{d}{dt} \overline{q}_{i0} = \m_1(\overline{q}_{i+1,0} - \overline{q}_{i0}) + \tilde{\lam } \mathds{1}\{q_{00} = 0\} \left[ \left(\overline{q}_{i-1,0}+\overline{q}_{i-1,1}\right)^{d_1} - \left(\overline{q}_{i0}+\overline{q}_{i-1,1}\right)^{d_1} \right]\\
\tfrac{d}{dt} \overline{q}_{i1} = \m_1(\overline{q}_{i+1,1} - \overline{q}_{i1}) + \tilde{\lam } \mathds{1}\{q_{00} = 0\} \left[ \left(\overline{q}_{i0}+\overline{q}_{i-1,1}\right)^{d_1} -\left(\overline{q}_{i0}+\overline{q}_{i1}\right)^{d_1}\right]
\end{dcases}
\end{equation}
with $\overline{q}_{00}+\overline{q}_{01}=1$.

Since the system operates under a preemptive priority discipline, the structure of the departure rate in each of the equations in~(\ref{eq:fluid_limit}) is clear. For instance, in order to change the proportion $\overline{q}_{11}$ due to a job completion, this job completion must take place at a server with configuration $(1,1)$. Exactly a fraction $q_{11} = \overline{q}_{11} - \overline{q}_{21}$ of the servers has this configuration and since these servers each work at rate $\m_1$ the total rate of change is given by $\m_1(\overline{q}_{21} - \overline{q}_{11})$.

 Let us illustrate the representation of the arrival term for the derivative of $\overline{q}_{11}$. Only an arrival of a type~I job at a server with configuration $(0,1)$ can contribute to the arrival term and the probability that this configuration is the smallest among the $d_1$ servers in the $d_1$-selection is given by $\left(\overline{q}_{10}+\overline{q}_{01}\right)^{d_1} -\left(\overline{q}_{10}+\overline{q}_{11}\right)^{d_1}$. Ties are broken according to the presence of a type~II job, in favor of having no type~II jobs.
 Moreover, there should be no idle servers because otherwise an arriving job would be allocated here as a type~II job.  
   Since type~I jobs arrive at a reduced rate~$\tilde{\lam}$, the total rate of change is given by  $\tilde{\lam } \mathds{1}\{q_{00} = 0\} [ \left(\overline{q}_{10}+\overline{q}_{01}\right)^{d_1} -\left(\overline{q}_{10}+\overline{q}_{11}\right)^{d_1}]$. \\
  
The expressions for the arrival terms in the derivatives of~(\ref{eq:fluid_limit}) and the reduced arrival rate~$\tilde{\lam}$ should be considered more carefully due to the discontinuity at $q_{00}=0$. We will give a sketch of the derivation of this fluid limit in Subsection~\ref{subsec:proofs_fl}. This derivation relies on the martingale method for point processes and Markovian queueing settings outlined by Pang \etal \cite{pang2007martingale} and Br\'{e}maud \cite{bremaud1981}. 

\begin{figure}[h]
\centering
\includegraphics{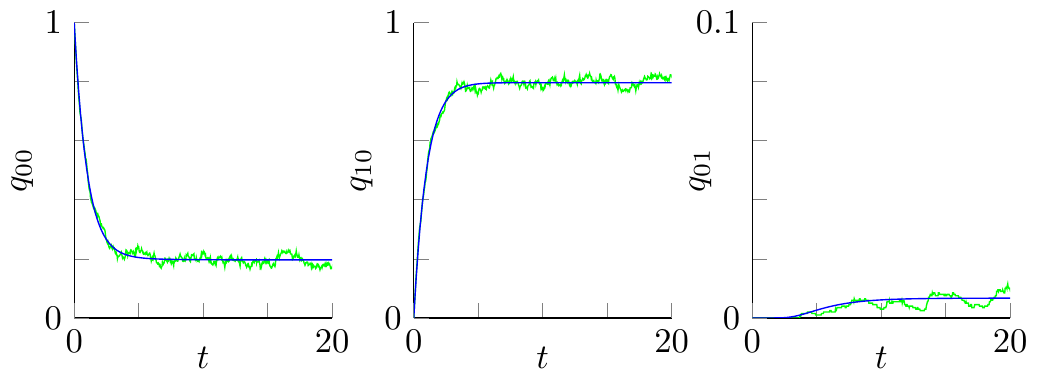}
\caption{A comparison between a simulated trajectory of the fluid-scaled stochastic process with $N=2000$ servers (green) and the numerical solution (blue) of the fluid limit (\ref{eq:fluid_limit}). The model parameters are given by $d_1 = 25$, $\lam = 0.8$, $\m_1 = 1$ and $\m_2 = 0.5$.}
\label{fig:sim_fluid_subplot_fav} 
\end{figure}

 The fluid-limit expression can be validated with simulations of the fluid-scaled stochastic process. Consider for instance Figure~\ref{fig:sim_fluid_subplot_fav} where the solution of the fluid limit~(\ref{eq:fluid_limit}) is presented together with a simulated trajectory of a reasonably large system. It can be observed that the simulated trajectory fluctuates closely around the numerical solution of the fluid limit, which supports the connection between the fluid limit and the behavior of the stochastic system in a many-server setting.

%%%%%%%%%%%%%%%%%%%%%%%%%%%%%%%%%%%%%%%%%%%%%%%%%%%%%%%%%%%%%%%%%%%%%%%%%%%%
%%%%%%%%%%%%%%%%%%%%%%%%%%%%%%%%%%%%%%%%%%%%%%%%%%%%%%%%%%%%%%%%%%%%%%%%%%%%%

\subsection{Fixed points} \label{subsec:fixedpoint}

To investigate the long-run behavior of the fluid limit~(\ref{eq:fluid_limit}), we are interested in its fixed points. It turns out that the mutual relationships between the model parameters $d_1$, $\lam$, $\m_1$ and $\m_2$ play a crucial role. 
 In the remainder of this section we investigate the setting where $\lam > \m_2$, in order to compare one of the fixed points with the fixed point of a JSQ($d_1$) policy with reduced load.  

\begin{theorem}\label{th:fixedpoint}\emph{(Fixed points).}
When $\lam > \m_2$ and $d_1\ge 2$, the system of differential equations~$(\ref{eq:fluid_limit})$ always has the following fixed point:
\begin{equation} \label{eq:fixedpoint}
%\left\{
\begin{cases}%{rcl}
\overline{q}_{i0}^*  = 0, &\\
\overline{q}_{i1}^* = \left(\frac{\lam-\m_2}{\m_1-\m_2}\right)^{\frac{d_1^i-1}{d_1-1}}, & i=0,1,2\dots.
\end{cases}
%\right.
\end{equation}
Furthermore, when $d_1 \ge d_1^*(\lam,\m_1,\m_2)$ precisely two more fixed points exist. These fixed points are such that $q_{00} + q_{01}+q_{10} = 1$ and $q_{00}>0$.  With $d_1^*\doteq d_1^*(\lam,\m_1,\m_2)$ the minimum selection size that satisfies
\begin{equation} \label{eq:fixedpoint_restrict}
\begin{array}{rcl}
d_1^* \lam \left( \frac{1}{\m_2} - \frac{1}{\m_1} \right) & > & 1, \\
\left(1- \frac{1}{d_1^*}\right) \frac{\m_1}{\lam} & > &  \left( d_1^*\lam \left( \frac{1}{\m_2} - \frac{1}{\m_1}\right) \right) ^{\frac{1}{d_1^*-1}}.
\end{array}
\end{equation}
\end{theorem}

The proof of this theorem can be found in Subsection~\ref{subsec:proofs_fl}.
 It can be observed that there always exists a sufficiently large $d_1$ value that satisfies both inequalities of condition~(\ref{eq:fixedpoint_restrict}) for given values of $\lam$, $\m_1$ and $\m_2$. This is trivial to see for the first inequality. The second inequality can be rewritten as
\begin{equation}
\left(1-\frac{1}{d_1}\right)\left(\frac{1}{d_1a}\right)^{\frac{1}{d_1-1}} > \frac{\lambda}{\m_1}
\end{equation}
with $a\doteq \lambda(\frac{1}{\m_2}-\frac{1}{\m_1})$.
The left hand side is increasing as function of $d_1$ with limit $1> \lambda/\m_1$. Table~\ref{tab:fp} gives the value of $d_1^*$, so that the condition~(\ref{eq:fixedpoint_restrict}) is satisfied for given model parameters and all $d_1\ge d_1^*$. It can be seen that the higher the load, the larger the size of the primary selections must be for multiple fixed points to persist. The additional fixed points have a strictly positive fraction of idle servers and it is intuitively clear that the number of servers where a job can be processed at rate $\m_1$ must grow with the load in order for these fixed points with a strictly positive fraction of idle servers to persist.\\

\begin{table}[h]
\centering
\begin{tabular}{c|cccccc}
$\lambda$ & 0.4 & 0.5 & 0.6 & 0.7 & 0.8 & 0.9 \\ \hline
$\mu_2=1/2$ & $/$ & $/$ & 5 & 9 & 18 & 46 \\
$\mu_2=1/3$ & 3 & 5 & 7 & 12 & 22 & 54
\end{tabular}
\caption{For given $\lam$, $\m_1 = 1$ and $\m_2$, the minimum value $d_1^*$ that satisfies condition~(\ref{eq:fixedpoint_restrict}) is listed.}
\label{tab:fp}
\end{table}

Let $\tilde{\lam}$ be as defined in~(\ref{eq:lam_tilde}) for the fixed point~(\ref{eq:fixedpoint}). 
 Then the long-term fraction of servers with at least $i$ jobs under a JSQ($d_1$) policy where each server works at rate $\m_1$ is given by
\begin{equation}
\overline{q}_i^* =  \left(\frac{\tilde{\lam}}{\m_1}\right)^{\frac{d_1^i-1}{d_1-1}} =  \left(\frac{\lam-\m_2}{\m_1-\m_2}\right)^{\frac{d_1^i-1}{d_1-1}},
\end{equation} 
for $i\ge 0$ \cite{mitzenmacher2001power}.
This shows a strong similarity with the fixed point~(\ref{eq:fixedpoint}) where two types of jobs are taken into account.
 Next we consider the case $d_1=1$. When $\lam > \m_2$ there still is a unique fixed point with $q_{00} = 0$, given by:
\begin{equation}
%\left\{
\begin{cases}%{rcl}
\overline{q}_{i0}^* = 0, &\\
\overline{q}_{i1}^* = \left(\frac{\lam-\m_2}{\m_1-\m_2}\right)^i, &i=0,1,2,\dots.
\end{cases}
%\right.
\end{equation}
This shows strong resemblance with the RA policy with load $\rho = \tilde{\lam}/\m_1$.
Allowing a primary selection of at least two servers leads to a super-exponential improvement compared to a primary selection of size one. On the other hand, there is no fixed point with $\lam \ge \m_2$ and $q_{00}>0$. Only if $\lam< \m_2$ we can show that there is a unique fixed point with $q_{00}>0$, namely
\begin{equation}
q_{00}^* = \frac{(\m_2-\lam)\m_1}{(\m_2-\lam)\m_1 + \lam\m_2}.
\end{equation}

\subsection{Further analysis} \label{subsec:analysis}
We will conduct a further analysis of the fluid limit~(\ref{eq:fluid_limit}) where we distinguish between $d_1$ sufficiently small and large compared to $d_1^*$ in the sense of conditions~(\ref{eq:fixedpoint_restrict}) in Theorem~\ref{th:fixedpoint}.

\subsubsection{Sufficiently small primary selections}

\begin{figure}[h]
\centering
\includegraphics{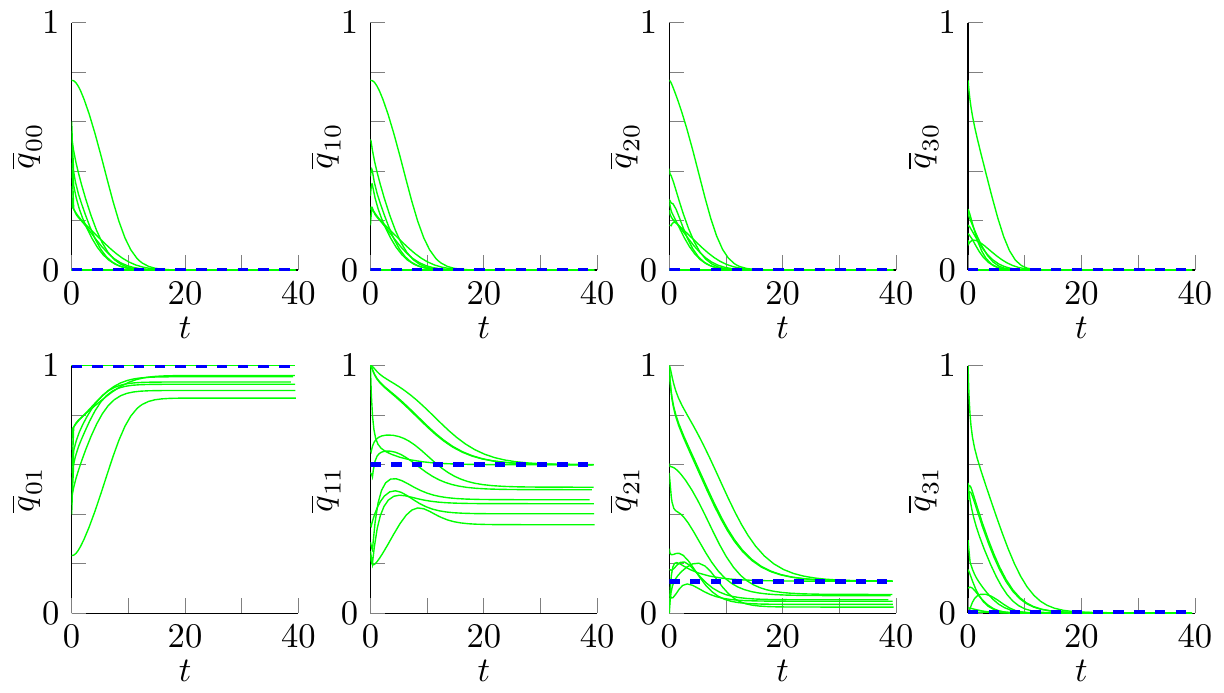} %[width=0.9\textwidth]
\caption{Ten solution trajectories (green) of the fluid limit (\ref{eq:fluid_limit}) are plotted for randomly generated initial states, showing convergence of the cumulative distributions to the fixed point (dashed, blue). The model parameters are given by $d_1 = 3$, $\lam  = 0.8$, $\m_1=1$ and $\m_2= 0.5$ and each of the servers has at most three type~I jobs in the initial states.}
\label{fig:globstab_unfav_subplot}
\end{figure}

When $d_1$ is sufficiently small in terms of the model parameters $\lam$, $\m_1$ and $\m_2$, i.e. $d_1 < d_1^*$, the fixed point~(\ref{eq:fixedpoint}) of the fluid limit~(\ref{eq:fluid_limit}) is unique. Numerical experiments suggest that this fixed point is a global attractor, i.e. the trajectories of the fluid limit will converge to this fixed point for every initial state of the system. As an example, we present Figure~\ref{fig:globstab_unfav_subplot} where the numerical solution of the fluid limit is visualized for ten randomly sampled initial configurations. We consider a system with the above-mentioned model parameters and a primary selection of size $d_1=3$.
As can be seen from the figure, all cumulative fractions $\left(\overline{q}_{i0}\right)_{i\ge0}$ tend to zero. Although some variability can be seen in the limiting behaviour of $(\overline{q}_{i1})_{i\ge0}$, the values are of the same order of magnitude as those of the theoretical fixed point. These deviations may be caused by numerical issues. There are two main reasons: 
(i) The fluid limit~(\ref{eq:fluid_limit}) is an infinite system of differential equations, so in order to solve the system numerically we need to truncate the system at some point. We choose to work with the variables up to $i=9$. 
(ii) The fluid limit~(\ref{eq:fluid_limit}) contains the indicator function $\mathds{1}\{q_{00} = 0\}$, while we use $\mathds{1}\{q_{00} < 10^{-15}\}$. This plays a role for initial conditions where a large fraction of the servers has an empty type~II queue, since the system still needs to make the transition to a state where every server has a type~II job.\\
 %%%%%%%%%%%%%%%%%%%%%%%%%%%%%%%%%%%%%%%%%%%%%%%%%%%%%%%
%%%%%%%%%%%%%%%%%%%%%%%%%%%%%%%%%%%%%%%%%%%%%%%%%%%%%%%%%%%%%%%%%%%%

\begin{figure}[h]
\centering
\includegraphics[width=9cm]{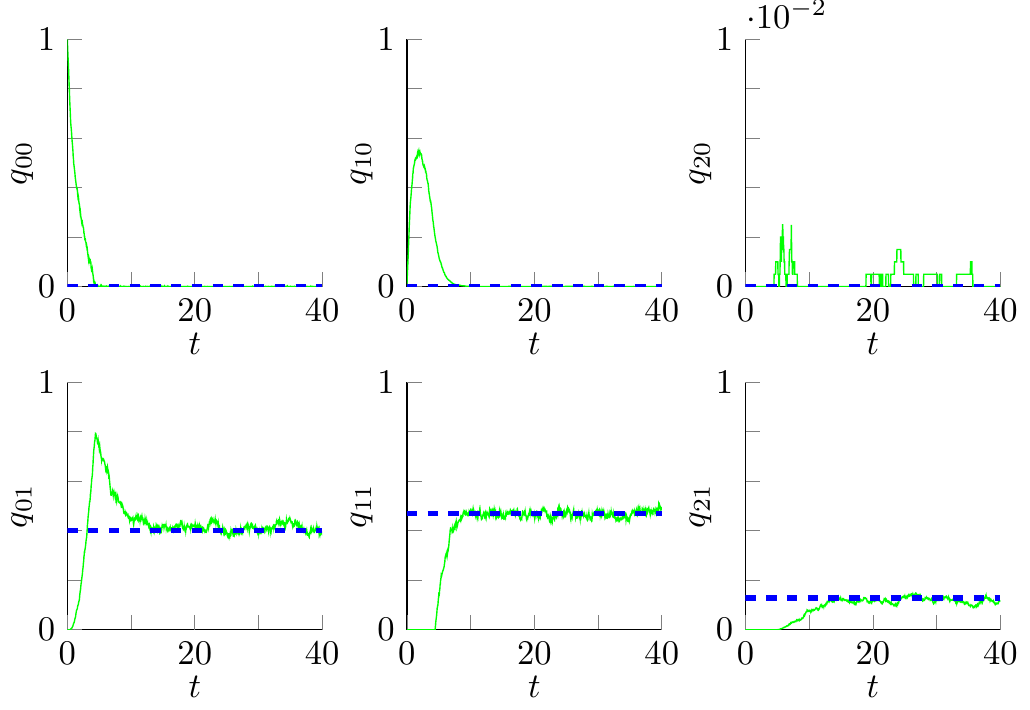}
\caption{A comparison between the long-term behavior of a simulated trajectory for $N= 2000$ servers (green) and the unique fixed point (blue) of the fluid limit (\ref{eq:fluid_limit}). The model parameters are given by $d_1 = 3$, $\lam = 0.8$, $\m_1 = 1$ and $\m_2 = 0.5$.}
\label{fig:sim_fluid_subplot_unfav}
\end{figure}

In the previous section we used the affinity coupling to show stochastic stability and the existence of an (unknown) stationary distribution for $\lam<\mu_1$. Assuming global stability of the unique fixed point, Theorem~1 by Bena\"{i}m and Le Boudec \cite{benaim2011mean} ensures that the large-$N$ limit of the stationary distribution will converge to the fixed point. 
Moreover, from simulations it can be observed that the trajectories converge to the unique fixed point of the fluid limit~(\ref{eq:fluid_limit}). As an example consider a system with above-mentioned model parameters. Figure~\ref{fig:sim_fluid_subplot_unfav} shows a simulated trajectory of the fluid-scaled variables for a system with $N=2000$ servers that is initially completely empty. It can be  seen that the trajectory converges to the fixed point $(q_{01},q_{11},q_{21},\dots) = (0.40,0.4704,0.1283,\dots)$, rounded at four decimals.\\

%%%%%%%%%%%%%%%%%%%%%%%%%%%%%%%%%%%%%%%%%%%%%%%%%%%%%%%%%%%%%%%%%%%%%%%%%%%%
%%%%%%%%%%%%%%%%%%%%%%%%%%%%%%%%%%%%%%%%%%%%%%%%%%%%%%%%%%%%%%%%%%%%%%%%%%%

 The asymptotic approximation for the mean stationary queue length, excluding the job in service, suggested by the fixed point is given by 
 \begin{equation}
 \E\left[Q_{ \textrm{CM($d_1$)}}\right] =  \suml_{i\ge 1} i q_{i,1}  =  \suml_{i\ge 1}  \left(\frac{\lam-\m_2}{\m_1-\m_2}\right)^{\frac{d_1^{i}-1}{d_1-1}}.
 \end{equation}
 Here CM($d_1$) refers to the combinatorial model with a primary server selection of size $d_1$.
 It is interesting to compare this with the asymptotic approximation for the mean queue length under a JSQ($d_1$) policy in the ordinary supermarket  model with arrival rate $\lam$ and service rate $\m_1$ \cite{mitzenmacher2001power,vvedenskaya1996queueing},
 \begin{equation}
 \E\left[Q_{\textrm{JSQ($d_1$)}} \right] =  \suml_{i\ge 1} i q_{i+1} =  \suml_{i\ge 1} \left(\frac{\lam}{\m_1}\right)^{\frac{d_1^{i+1}-1}{d_1-1}},
 \end{equation}
 and the exact mean queue length under the RA policy,
 \begin{equation}
 \E\left[Q_{\mathrm{RA}}\right]  = \left(1-\frac{\lam}{\m_1}\right) \suml_{i\ge 1} i \left(\frac{\lam}{\m_1}\right) ^{i+1}  = \frac{(\lam/\m_1)^2}{1-\lam/\m_1}.
 \end{equation}
 
\begin{figure}[h]
\centering
\includegraphics[width=10cm]{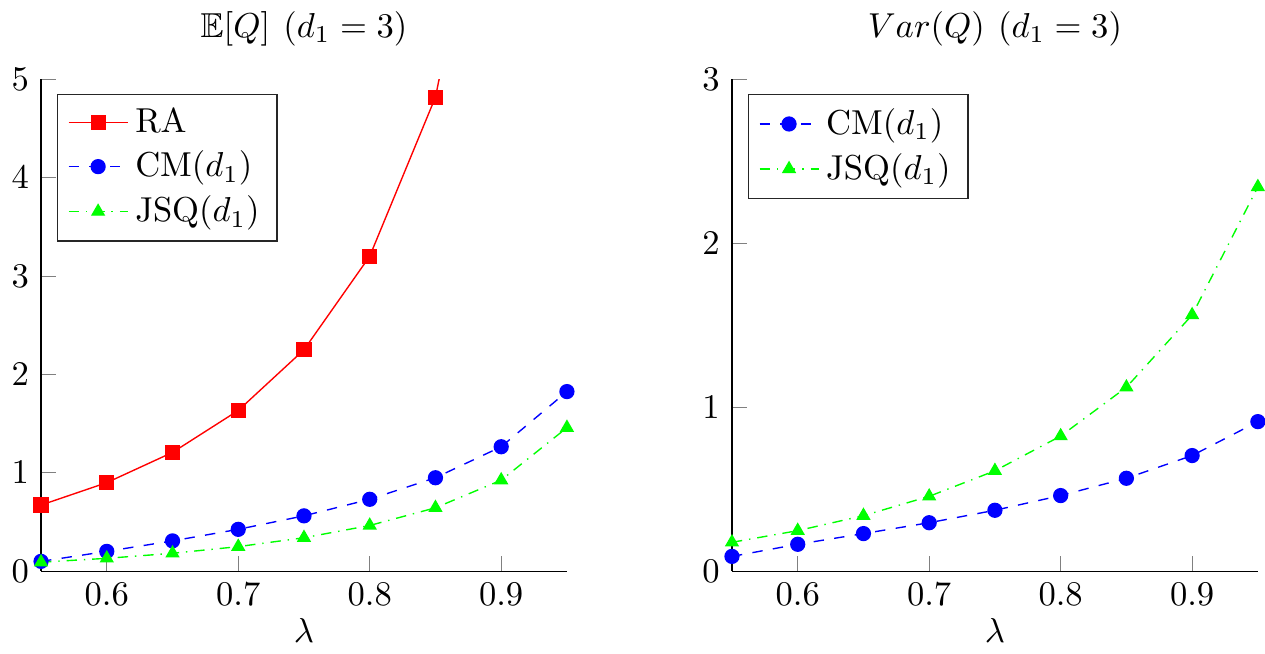}
\caption{Comparison of mean and variance of the queue length as a function of $\lam$ between three different models on fluid level for $d_1 =3$, $\m_1 = 1$ and $\m_2 = 0.5$.}
\label{fig:compare_model} 
\end{figure}

 Figure~\ref{fig:compare_model} presents a comparison of the number of waiting jobs as a function of $\lam$, with $d_1=3$, $\m_1=1$ and $\m_2=0.5$. It is known that the mean queue length for the RA policy tends to infinity when the offered traffic grows to one. We see that the mean queue length in the combinatorial model is slightly larger than for the JSQ($d_1$) policy. On the other hand the variance of the queue length in the JSQ($d_1$) model is almost twice as large compared to the combinatorial model. We conclude that the combinatorial model still performs well from a queue length perspective, even though each server has a type~II job and possibly multiple type~I jobs in its queue.

%%%%%%%%%%%%%%%%%%%%%%%%%%%%%%%%%%%%%%%%%%%%%%%%%%%%%%%
%%%%%%%%%%%%%%%%%%%%%%%%%%%%%%%%%%%%%%%%%%%%%%%%%%%%%%%%%%%%%%%%%%%%

 From the fixed point expression, it is not immediately visible that type~II jobs finish their service since the fraction $q_{00}$ is zero. However, an idle server will be filled instantly with an arriving job. A total fraction 
\begin{equation}
\frac{\lam-\tilde{\lam}}{\lam} = \frac{\m_2}{\lam}\frac{\m_1-\lam}{\m_1-\m_2}
\end{equation}
of the arriving jobs undergo this `immediate switch': they are allocated as a type~II job to a server that just emptied its queue. This fraction decreases in $\lam$ and for example for the above-mentioned model parameters, this leads to a fraction of $1/4$. 
Furthermore, the type~II jobs will leave the system at the same rate $\lam -\tilde{\lam}$ as they enter the system, since we study the system in equilibrium. Moreover, due to Little's law we know that the expected waiting time of an arbitrary job is finite. Let $W$ denote the waiting time, then
$
\E\left[Q_{\textrm{CM($d_1$)}} \right]  = \lam \E[W].
$
Since the expected queue length under our affinity-scheduling policy is finite, this results in a finite expected waiting time for an arbitrary job, so also for the type~II jobs.\\

\begin{figure}[h]
\centering
\includegraphics[width=8cm]{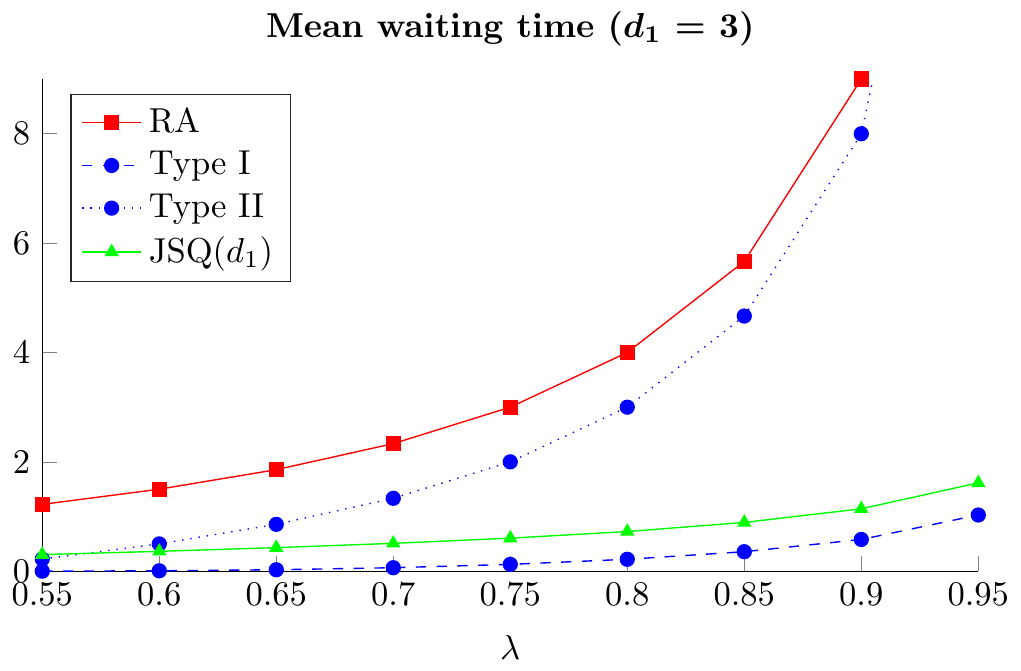} %[width=0.9\textwidth]
\caption{A comparison of the mean waiting times of the type~I and type~II jobs (blue) with the mean waiting times under the RA policy (red) or a JSQ($d_1$) policy (green) as function of $\lam$ for $d_1 =3$, $\m_1 = 1$ and $\m_2 = 0.5$.}
\label{fig:compare_model_WT} 
\end{figure}

Since each server operates under a preemptive scheduling policy, we can calculate the average waiting time of a type~I job using Little's law. Let $Q_{\mathrm{I}}$ denote the number of type~I jobs at a server. Then 
\begin{equation}
\E[Q_{\mathrm{I}}] = \suml_{i\ge 1} i q_{i+1,1}  =  \suml_{i\ge 1} \left(\frac{\lam-\m_2}{\m_1-\m_2}\right)^{\frac{d_1^{i+1}-1}{d_1-1}} .
\end{equation}
Furthermore, the reduced arrival rate $\tilde{\lam}$ gives the arrival rate of type~I jobs on fluid level. If $W_\mathrm{I}$ represents the waiting time of a type~I job, then due to Little's law
$
\E[Q_{\mathrm{I}}] = \tilde{\lam}\E[W_{\mathrm{I}}].
$
Let $Q_{\mathrm{II}}$ and $W_{\mathrm{II}}$ have the same interpretation as above but for the type~II jobs. We condition on the type of job to obtain
\begin{equation}
\E[W] = \frac{\tilde{\lam}}{\lam} \E[W_{\mathrm{I}}] + \frac{\lam - \tilde{\lam}}{\lam} \E[W_{\mathrm{II}}].
\end{equation}
Because of Little's law this
results in 
\begin{equation}
\E[W_{\mathrm{II}}] = \frac{1}{\lam-\tilde{\lam}} \left( \E\left[Q_{\mathrm{CM\left(\right.}d_1\mathrm{\left.\right)}} \right]  - \E[Q_{\mathrm{I}}]\right) = \frac{1}{\lam-\tilde{\lam}} \frac{\lam-\m_2}{\m_1-\m_2}.
\end{equation}
We can also immediately apply Little's law to the type~II jobs. We know that they arrive at rate $\lam-\tilde{\lam}$ and the mean waiting queue length is by definition given by
\begin{equation}
\E[Q_{\mathrm{II}}] = \suml_{i \ge 1} q_{i1} = 1- q_{01} = \frac{\lam-\m_2}{\m_1-\m_2}.
\end{equation}
In Figure~\ref{fig:compare_model_WT} we compare the mean waiting time of a type~I or type~II job with the mean waiting time under the RA or the JSQ($d_1$) policy. The mean waiting time of type~II jobs is fairly high, but still lower than the waiting time under the RA policy. We also observe that the mean waiting time of type~I jobs is significantly smaller than under a JSQ($d_1$) policy. We conclude that our allocation strategy leads to a reduction in the mean waiting time for a large group of arriving jobs at the expense of some other jobs that encounter longer waiting times. The uniqueness of the fixed point allows us to analyze the asymptotic stationary distribution of the model, on the other hand we observe that the value of the size of the server selection $d_1$ is too small to achieve a zero waiting time for an arriving job.

%%%%%%%%%%%%%%%%%%%%%%%%%%%%%%%%%%%%%%%%%%%%%%%%%%%%%%%
%%%%%%%%%%%%%%%%%%%%%%%%%%%%%%%%%%%%%%%%%%%%%%%%%%%%%%%%%%%%%%%%%%%%

\subsubsection{Sufficiently large primary selections}

Assume that the primary selection has a sufficiently large size $d_1$ for given model parameters in terms of the conditions~(\ref{eq:fixedpoint_restrict}), i.e. $d_1\ge d_1^*$. From Theorem \ref{th:fixedpoint} we know that, next to the closed form fixed point (\ref{eq:fixedpoint}), there are two additional fixed points with $q_{00}+q_{01}+q_{01}=1$. We prove the following theorem using the indirect Lyapunov method.

\begin{theorem}\label{th:fixed_unst}\emph{(Local (in)stability).}
Of the two additional fixed points mentioned in Theorem~\ref{th:fixedpoint} with $q_{00}+q_{01}+q_{01}=1$ when $d_1\ge d_1^*$, one is locally stable and the other one is unstable.
\end{theorem}

\begin{figure}[h]
\centering
\includegraphics[width=8cm]{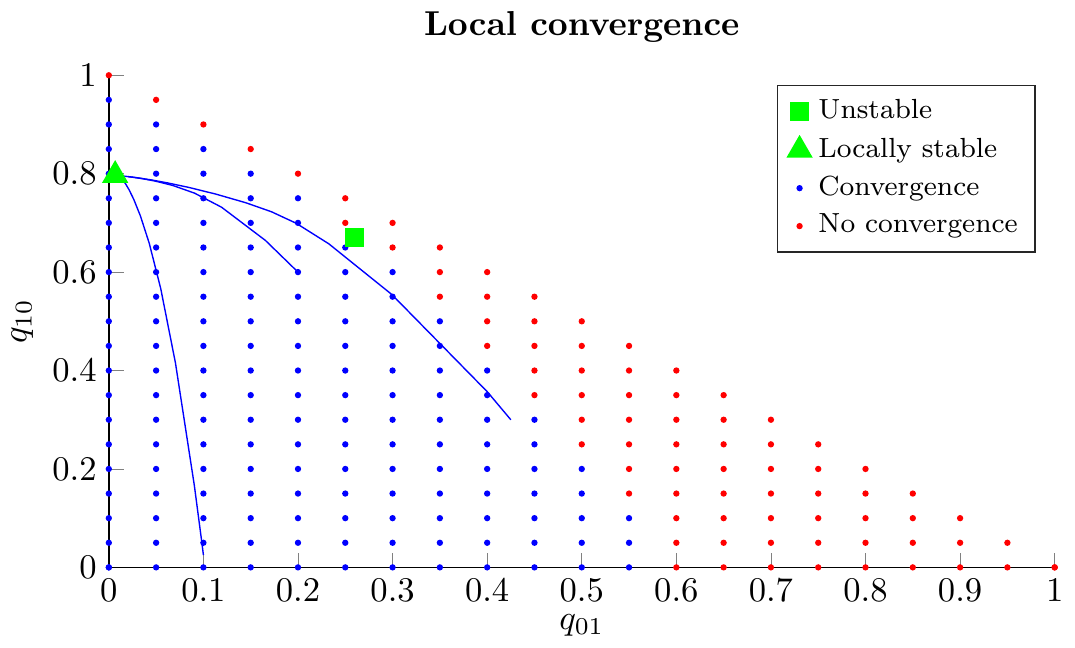} 
\caption{An overview of the initial states with $q_{00}+q_{01}+q_{10} = 1$, if the corresponding trajectories converge to the locally stable fixed point (green triangle), the initial states are indicated in blue otherwise they are indicated in red. We consider a system with $d_1=25$, $\lam=0.8$, $\m_1=1$ and $\m_2=0.5$.}
\label{fig:localstab} 
\end{figure}

The proof of Theorem~\ref{th:fixed_unst} is given in Subsection~\ref{subsec:proofs_fl}.
In the remainder of this subsection, we will provide a numerical illustration, where we consider a system with $\lam=0.8$, $\m_1=1$, $\m_2=0.5$ and $d_1=25$ throughout. We observed similar qualitative behavior  across many different scenarios, but only present results for those parameter values because of space constraints.
To get a better notion of the local stability we present Figure~\ref{fig:localstab}. For several initial values such that $q_{00}+ q_{01}+q_{10} = 1$, the system of differential equations~(\ref{eq:fluid_limit}) is solved numerically. %We consider a system where $d_1$, $\lam$, $\m_1$ and $\m_2$ are 25, 0.8, 1, 0.5, respectively.
 All trajectories with initial states indicated in blue will converge to the locally stable fixed point from the previous theorem and a few of these trajectories are also visualized. All other initial states, indicated in red, will not converge to this locally stable fixed point. We see that these states have a large fraction of servers with a type~II job present and a small fraction of idle servers, since there is a smaller probability to select an idle server in the $d_1$-selection. So jobs will have a longer mean service time as a type~II job and jobs will start to accumulate.
 
%%%%%%%%%%%%%%%%%%%%%%%%%%%%%%%%%%%%%%%%%%%%%%%%%%%%%%%
%%%%%%%%%%%%%%%%%%%%%%%%%%%%%%%%%%%%%%%%%%%%%%%%%%%%%%%%%%%%%%%%%%%%
 
\begin{figure}[h]
\centering
\includegraphics[width=8cm]{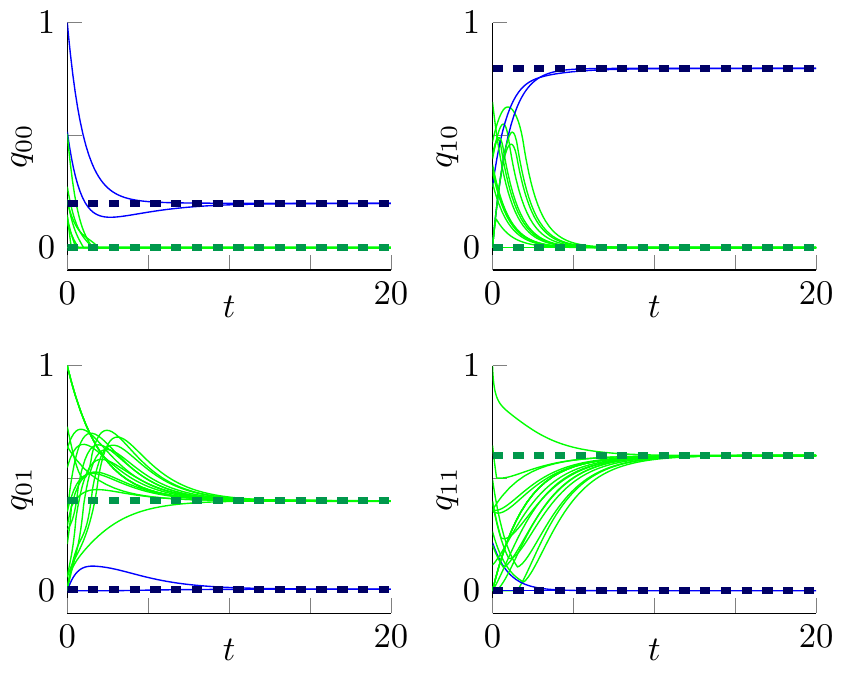} 
\caption{Trajectories for 20 randomly sampled initial states are visualized, with initial points such that  $q_{00}+q_{01}+q_{10}+q_{11} = 1$. Trajectories in blue converge to the no-queueing fixed point while green trajectories will converge to the queueing fixed point. We consider a system where $d_1=25$, $\lam=0.8$, $\m_1=1$ and $\m_2=0.5$, and both fixed points are indicated by dashed lines.}
\label{fig:two_attract_subplot} 
\end{figure}

 In total this gives rise to two locally stable fixed points: the closed-form fixed point~(\ref{eq:fixedpoint}) where each server has a type~II job and possibly multiple type~I jobs and the fixed point from Theorem~\ref{th:fixed_unst} where at most one job is present at each server. In the remainder of this section we will refer to these fixed points as the \textit{queueing fixed point} and \textit{no-queueing fixed point}, respectively. 
 We do not formally prove this statement but we will illustrate it with a representative example.  
For a system with the above-mentioned parameters, the two fixed points under consideration (non-cumulative fractions) are given by:
\begin{equation} \label{eq:fixedpoint_filledin}
\begin{array}{rcl}
(q_{00},q_{01},q_{10}) & = & (0.1966,0.0067,0.7967)\\
(q_{01},q_{11}) & = & (0.4,0.6).
\end{array}
\end{equation}
Both fixed points are indicated with dashed lines in Figure~\ref{fig:two_attract_subplot}, in dark blue and dark green, respectively. Furthermore, the graphs contain 20 trajectories starting from randomly sampled initial configurations with $q_{00}+q_{01}+q_{10}+q_{11} = 1$, all these trajectories converge to one of the two fixed points. This implies that the convergence area presented in Figure~\ref{fig:localstab} to the single-user fixed point will in fact be larger. As can be seen, most of the trajectories will converge to the type-II fixed point. This phenomenon will be even more apparent if we allow initial states with more than two jobs.

%%%%%%%%%%%%%%%%%%%%%%%%%%%%%%%%%%%%%%%%%%%%%%%%%%%%%%%%%
%%%%%%%%%%%%%%%%%%%%%%%%%%%%%%%%%%%%%%%%%%%%%%%%%%%

The literature often describes systems with a unique global attractor as a fixed point of the fluid limit so that there is a direct connection between the stationary distribution in a many-server setting and this fixed point. However, the non-uniqueness of the fixed points does not imply that these two concepts are completely uncorrelated. For instance, Figure~\ref{fig:sim_fluid_subplot_fav} presents a comparison between the numerical solution of the fluid limit and a simulation with $N=2000$ servers with the above-mentioned model parameters. The system is initially empty and the simulated trajectory seems to converge to the no-queueing fixed point. 
We can presents a similar figure, where in the initial configuration each server has one type~II job, in which case both the numerical solution and the simulation seem to tend to the queueing fixed point.\\

 However, the stochastic process with a finite number of servers is an irreducible Markov process which implies that any state can be reached as long as the process is observed long enough and a unique equilibrium distribution must exist. 
 Nevertheless, it can be observed that the process spends a long 
 Nevertheless, it can be observed that the residence time near each of the locally stable fixed points, which increases with $N$, is long before the process makes the transition to the other locally stable fixed point. Gibbens \etal \cite{gibbens1990bistability} describe this concept of switching between multiple modes by `tunneling'.
 We can call these locally stable points the `quasi-stationary' distributions of the stochastic process as in \cite{zachary2011loss}.
 
 Examples of models with multiple local fixed points in loss and communication networks can be found in \cite{bean1997dynamic,gibbens1990bistability,zachary2011loss}. More recent work by Martirosyan and Robert~\cite{martirosyan2018equilibrium} considers an allocation strategy closely related to the affinity-scheduling policy in a loss network setting, i.e. jobs can be redirected to distant servers with a penalty or can be omitted if none of the servers has enough spare capacity. Also in this setting, a fluid-limit analysis reveals multiple locally stable fixed points.

%%%%%%%%%%%%%%%%%%%%%%%%%%%%%%%%%%%%%%%%%%%%%%%%%%%%%%%%%%%%%%%%%%%%%%%%%%%%%%%%%%%%%%%%%%%%%%%%%%%%%%%%%%%%%%%%%%%%%%%%%%%%%%%%%%%%%%%%%%%%%%%%%%%%%%%%%% 
\section{Proofs} \label{sec:proofs}

\subsection{Proof of Lemma~\ref{lem:generalcoupling}: Affinity coupling} \label{subsec:proof_stoch_dom}
Since the system configurations between two consecutive events remain unchanged, we will condition on the discrete event times and use forward induction.\\
Assume that~(\ref{eq:lem_maj}) holds up to the time of the $(k-1)$-th event, we will argue that the majorization property still holds at time $t_k$ of the $k$-th event by making a distinction between arrival and departure epochs. But first we need a formal way to express the effect of these events in terms of $(\overline{Q}_{i}^{\mathrm{aff}}(t) )_{i\ge1}$ and $(\overline{Q}_{i}^{\mathrm{ref}}(t) )_{i\ge1}$. For instance, let $n$ be the server position selected for a departure. Due to the ordering we know that there are at least $N-n+1$ servers with the same number of jobs or more in their queues as the server at position $n$. It might also be possible that the server at position $n-1$ has the same number of jobs as the server at position $n$, there is notable difference in in terms of the variables $\overline{Q}_{i}$ whether a removal takes place at position $n-1$ or at position $n$. Instead of removing from the server at position $n$ and reordering the servers before computing $(\overline{Q}_{i}^{\mathrm{aff}}(t) )_{i\ge1}$ and $(\overline{Q}_{i}^{\mathrm{ref}}(t) )_{i\ge1}$, we can also immediately compute these quantities. The difference is subtle and valid because the proof does not rely on the present type~II jobs or on the actual servers but only on their relative positions. Therefore we define two intermediate quantities: 

\begin{equation}\label{eq:defprop1}
\begin{array}{rcl}
I_{\mathrm{aff}}(n) & := & \max\{j: \overline{Q}_{j}^{\mathrm{aff}}\ge N-n+1\}, \\
I_{\mathrm{ref}}(n) & := & \max\{j: \overline{Q}_{j}^{\mathrm{ref}}\ge N-n+1\}.
\end{array}
\end{equation}
For instance in the original system in Figure~\ref{fig:sample_path}, $I_{\mathrm{aff}}(n)$ is given by 3.
%This concept is visualized in Figure~\ref{fig:fig_defn}. 
Furthermore, only one job will be added or removed at a discrete time event. A new event at time $t_k$ could only violate~(\ref{eq:lem_maj}) if at time $t_{k-1}$ (\ref{eq:lem_maj})~holds with equality, i.e. 
\begin{equation} \label{eq:graph_majo_eq}
\suml_{i = m}^{ \infty}\overline{Q}_{i}^{\mathrm{aff}}(t_{k-1}) = \suml_{i = m}^{ \infty}\overline{Q}_{i}^{\mathrm{ref}}(t_{k-1})
\end{equation}
with $m \ge 1$. Therefore we only focus on this setting in the induction step.% while deriving the induction step.\\

\textit{Arrival.} At time $t_k$ an arrival occurs and first position $n$ is selected, the updated reference system looks as follows:
\begin{equation}
\overline{Q}_j^{\mathrm{ref}}(t_k) = 
\left\{
\begin{array}{ll}
\overline{Q}_j^{\mathrm{ref}}(t_k^-)+1, & \text{if } j = I_{\mathrm{ref}}(n)+1\\
\overline{Q}_j^{\mathrm{ref}}(t_k^-), & \text{otherwise.} \\
\end{array}
\right.
\end{equation}
If the newly arrived job is allocated as a type~II job in the original system or no arrival takes place due to the coupling, (\ref{eq:lem_maj}) is trivially satisfied. We consider the setting where the job is allocated as a type~I job to a server at position $n_{\mathrm{aff}}$ which is at most $n$, such that
\begin{equation}
\overline{Q}_{j}^{\mathrm{aff}}(t_k) = 
\left\{
\begin{array}{ll}
\overline{Q}_{j}^{\mathrm{aff}}(t_k^-)+1, & \text{if } j = I_{\mathrm{aff}}(n_{\mathrm{aff}})+1\\
\overline{Q}_{j}^{\mathrm{aff}}(t_k^-), & \text{otherwise.} \\
\end{array}
\right.
\end{equation}
Moreover, the left hand side of~(\ref{eq:lem_maj}) remains unchanged if $m > I_{\mathrm{aff}}(n_{\mathrm{aff}})+1$ so that the order in~(\ref{eq:lem_maj}) is preserved. Now, fix $m \le I_{\mathrm{aff}}(n_{\mathrm{aff}})+1$, if we now show that also $I_{\mathrm{ref}}(n) \ge m-1$, then~(\ref{eq:lem_maj}) remains valid since both sides are raised by one. We use~(\ref{eq:graph_majo_eq}) and the induction hypothesis for $m-1$ at time $t_k^-$ to obtain
\begin{equation}
\begin{array}{rcccl}
\overline{Q}_{m-1}^{\mathrm{aff}}(t_k^-)& =& \suml_{i = m-1}^{ \infty}\overline{Q}_{i}^{\mathrm{aff}}(t^-_k)  - \suml_{i = m}^{ \infty}\overline{Q}_{i}^{\mathrm{aff}}(t^-_k) & &\\
& \le& \suml_{i = m-1}^{ \infty}\overline{Q}_{i}^{\mathrm{ref}}(t^-_k)  - \suml_{i = m}^{ \infty}\overline{Q}_{i}^{\mathrm{ref}}(t^-_k) & = & \overline{Q}_{m-1}^{\mathrm{ref}}(t_k^-).
\end{array}
\end{equation}
Then it follows that $I_{\mathrm{aff}}(n_{\mathrm{aff}}) \ge m-1$ implies $I_{\mathrm{ref}}(n) \ge m -1$
which concludes the derivation if the event at time $t_k$ is an arrival.

\textit{Departure.} If at time $t_k$ a departure will take place, one of the following four scenarios will occur.
\begin{enumerate}
\item There is a job completion of a type~I job in the original system and of a job in the reference system.
\item There is only a departure at the jobs of the reference system.
\item There is only a departure of a type~I job in the original system.
\item There is no departure at the type~I jobs of the original system or the jobs of the reference system.
\end{enumerate}
It is clear that we only need to investigate the first two scenarios.\\
\textit{Scenario 1.} Let $n \in W$ be the position of the servers in both the original and the reference system from which a job will be removed. The updated systems will look as follows,
\begin{equation}
\begin{array}{c}
\overline{Q}_{j}^{\mathrm{aff}}(t_k) = 
\left\{
\begin{array}{ll}
\overline{Q}_{j}^{\mathrm{aff}}(t_k^-)-1, & \text{if } j = I_{\mathrm{aff}}(n)\\
\overline{Q}_{j}^{\mathrm{aff}}(t_k^-), & \text{otherwise}, \\
\end{array}
\right.

\\

\overline{Q}_{j}^{\mathrm{ref}}(t_k) = 
\left\{
\begin{array}{ll}
\overline{Q}_{j}^{\mathrm{ref}}(t_k^-)-1, & \text{if } j = I_{\mathrm{ref}}(n)\\
\overline{Q}_{j}^{\mathrm{ref}}(t_k^-), & \text{otherwise}. \\
\end{array}
\right.
\end{array}
\end{equation}
We will focus on $m\le I_{\mathrm{ref}}(n)$, since for $m > I_{\mathrm{ref}}(n)$ (\ref{eq:lem_maj})~remains trivially valid. A similar argument as above will be used to show that $ I_{\mathrm{aff}}(n) \ge m$, so that both sides will be lowered by one compared to the event time $t_{k-1}$. We use~(\ref{eq:graph_majo_eq}) and the induction hypothesis for $m+1$ at time $t_k^-$ to obtain 
$\overline{Q}_{m}^{\mathrm{aff}}(t_k^-) \ge  \overline{Q}_{m}^{\mathrm{ref}}(t_k^-)$.
Then it follows that $I_{\mathrm{ref}}(n) \ge m$ implies $I_{\mathrm{aff}}(n) \ge m$
which concludes the proof of scenario~1.\\
\textit{Scenario 2.} Let  $n \in W_{\mathrm{ref}} \setminus W$ be the position where a job leaves the reference system, then for all $j$
\begin{equation}
\begin{array}{c}
\overline{Q}_{j}^{\mathrm{aff}}(t_k) = \overline{Q}_{j}^{\mathrm{aff}}(t_k^-),  \\
\\
\overline{Q}_{j}^{\mathrm{ref}}(t_k) = 
\left\{
\begin{array}{ll}
\overline{Q}_{j}^{\mathrm{ref}}(t_k^-)-1, & \text{if } j = I_{\mathrm{ref}}(n)\\
\overline{Q}_{j}^{\mathrm{ref}}(t_k^-), & \text{otherwise}. \\
\end{array}
\right.
\end{array}
\end{equation}
Again we focus on $m \le  I_{\mathrm{ref}}(n)$.
 Fix $m$, we will show by contradiction that~(\ref{eq:graph_majo_eq}) cannot occur so that~(\ref{eq:lem_maj}) is preserved at time $t_k$ since the right hand side can be lowered by at most one. Assuming that~(\ref{eq:graph_majo_eq}) does hold and using the induction hypothesis on $m+1$, we conclude that
 $\overline{Q}_{m}^{\mathrm{aff}}(t_k^-)\ge \overline{Q}_{m}^{\mathrm{ref}}(t^-_k)$.
 Now,
 \begin{equation}
 \overline{Q}_{m}^{\mathrm{aff}}(t^-_k) \ge \overline{Q}_{m}^{\mathrm{ref}}(t^-_k) \ge \overline{Q}_{I_{\mathrm{ref}}(n)}^{\mathrm{ref}}(t^-_k) \ge N-n+1 \ge |W|+1,
\end{equation}  
since $N-|W_{\mathrm{ref}}| < n\le N-|W|$. This implies that $\overline{Q}_{m}^{\mathrm{aff}}(t_k^-)  > |W|$, however there are only $|W| = |W_{\mathrm{aff}}|$ servers working on a type~I job in the original system. This leads to a contradiction and concludes the proof of Lemma~\ref{lem:generalcoupling}.

%The result of Lemma~\ref{lem:generalcoupling} is similar to the result of Proposition~3.1 in Mukherjee \etal \cite{mukherjee2016universality}. Both the original and the reference system are here operating on $N$~exchangeable servers, but the argumentation will also hold for more general modeling frameworks as described above.

%%%%%%%%%%%%%%%%%%%%%%%%%%%%%%%%%%%%%%%%%%%%%%%%%%%%%%%%%%%%

\subsection{Proofs: Fluid limit and fixed point analysis} \label{subsec:proofs_fl}

\subsubsection{Derivation fluid limit~$(\ref{eq:fluid_limit})$} \label{subsec:proof_fluidlimit}

First, consider the stochastic process with $N$ servers and its corresponding flow conservation equations. Next, the martingale methods as outlined by Pang \etal \cite{pang2007martingale} and Br\'{e}maud \cite{bremaud1981} are applied and the limit of $N$ to infinity of the fluid scaled process is studied. Then~(\ref{eq:fluid_limit}) is obtained from the resulting system of integral equations.\\

\textit{Step 1: flow conservation equations.} Let $p_{ij}^N(q,t)$ be the probability that an arriving job at time $t$ is allocated to a server with $i$ type~I jobs and $j$ type~II jobs as a type~$q$ job, with $q \in \{\mathrm{I},\mathrm{II}\}$. As before, we will omit the time dependence $t$ to ease the notation.

Only to an idle server we can allocate a job as a type~I or type~II job; allocations to servers with a higher configuration will always take place as a type~I job. The corresponding transition probabilities are given by
\begin{equation}
p_{00}^N(\mathrm{I}) = 1- \left(1-\frac{Q_{00}^N}{N}\right)^{d_1},
\end{equation} 
the probability that an idle server is present in the $d_1$-selection, and
\begin{equation}
p_{00}^N(\mathrm{II})  =  \mathds{1}\{Q_{00}^N >0 \} \left(1-\frac{Q_{00}^N}{N}\right)^{d_1},
\end{equation}
the probability that the $d_1$-selection does not contain an idle server while they are present. As mentioned in the model description, the $d_2$-selection contains all servers that are not in the $d_1$-selection. Hence the indicator  function $\mathds{1}\{Q_{00}^N >0 \}$ emerges in the probabilities.

An arriving job will be allocated as a type~I job to a server with configuration $(i,0)$, with $i\ge 1$, if the minimum configuration in the $d_1$-selection is given by $(i,0)$ and when there are no completely idle servers that can be included in the $d_2$-selection. The corresponding probability is given by the probability that the $d_1$-selection contains only servers with at least $i$ type~$I$ jobs minus the probability that all $d_1$ servers have a configuration strictly higher than $(i,0)$. Thus, for $i\ge 1$,
\begin{equation}
p_{i0}^N(\mathrm{I}) = \mathds{1}\{Q^N_{00} = 0\} \left[ \left( \suml_{k\ge i} \left[ \frac{Q^N_{k0}}{N} + \frac{Q^N_{k1}}{N} \right]\right)^{d_1} - \left( \frac{Q^N_{i1}}{N} + \suml_{k\ge i+1} \left[ \frac{Q^N_{k0}}{N} + \frac{Q^N_{k1}}{N}\right]\right)^{d_1 }\right].
\end{equation}
In a similar way, we obtain $p_{i1}^N(I)$, for $i\ge 0$:
\begin{equation}
p_{i1}^N(\mathrm{I}) = \mathds{1}\{Q^N_{00} = 0\} 
\left[ 
\left(  \frac{Q^N_{i1}}{N} +\suml_{k\ge i+1} \left[ \frac{Q^N_{k0}}{N} + \frac{Q^N_{k1}}{N}\right]\right)^{d_1}
 - \left(  \suml_{k\ge i+1} \left[ \frac{Q^N_{k0}}{N} + \frac{Q^N_{k1}}{N}\right]\right)^{d_1}
 \right].
\end{equation}
Once these probabilities are set, the flow conservation equations can be constructed. The randomness in the stochastic model is caused by Poisson arrivals and exponentially distributed service times, so that the number of arrivals and service completions can be counted using Poisson processes with appropriately chosen rates. Define a set of independent Poisson processes with rate 1. Let $P_{A_{00,q}}$ denote the Poisson counting process for the number of arriving type~$q$ jobs at servers with configuration $(0,0)$, and $P_{A_{ij}}~i+j \ge 1$ reflects the arriving jobs at servers with configuration $(i,j)$. Similarly, define the counting process of the service completions $P_{S_{ij}},~i+j \ge 1$.
Furthermore if $i \ge1 $, the number of servers at time $t$ with at least $i$ type~I jobs and  exactly $j$ type~II jobs depends on its initial state $(\overline{Q}^N_{ij}(0))$, the number of service completions of jobs at servers with configuration $(i,j)$ and the number of arrivals at servers in configuration $(i-1,j)$ within the time interval $[0,t)$. We obtain the following flow conservation equations for the stochastic model $(\overline{Q}^N_{ij})_{i,j}$ with $N$ servers and total arrival rate $\lam N$. Let $i\ge 2$:
\begin{equation}
\begin{array}{rcl}
\overline{Q}^N_{00}(t) & = & \overline{Q}^N_{00}(0) + P_{S_{01}}\left(\m_2 \intt Q_{01}^N(s)\,ds\right) - P_{A_{00}}\left(\lam N \intt \left[p_{00}^N(\mathrm{I},s) +  p_{00}^N(\mathrm{II},s)\right]\,ds\right)\\

\overline{Q}^N_{01}(t) & = & \overline{Q}^N_{01}(0) -  P_{S_{01}}\left(\m_2 \intt Q_{01}^N(s)\,ds\right) + P_{A_{00,\mathrm{II}}}\left(\lam N \intt  p_{00}^N(\mathrm{II},s)\,ds\right)\\

\overline{Q}^N_{10}(t) & = & \overline{Q}^N_{10}(0) -  P_{S_{10}}\left(\m_1 \intt Q_{10}^N(s)\,ds\right) + P_{A_{00,\mathrm{I}}}\left(\lam N \intt  p_{00}^N(\mathrm{I},s)\,ds\right)\\

\overline{Q}^N_{11}(t) & = & \overline{Q}^N_{11}(0) -  P_{S_{11}}\left(\m_1 \intt Q_{11}^N(s)\,ds\right) + P_{A_{01}}\left(\lam N \intt  p_{01}^N(\mathrm{I},s)\,ds\right)\\

\overline{Q}^N_{i0}(t) & = & \overline{Q}^N_{i0}(0) -  P_{S_{i0}}\left(\m_1 \intt Q_{i0}^N(s)\,ds\right) + P_{A_{i-1,0}}\left(\lam N \intt  p_{i-1,0}^N(\mathrm{I},s)\,ds\right)\\

\overline{Q}^N_{i1}(t) & = & \overline{Q}^N_{i1}(0) -  P_{S_{i1}}\left(\m_1 \intt Q_{i1}^N(s)\,ds\right) + P_{A_{i-1,1}}\left(\lam N \intt  p_{i-1,1}^N(\mathrm{I},s)\,ds\right).\\

\end{array}
\end{equation}
Due to the Poisson split property we define $P_{A_{00}}$ as the sum of the two processes $P_{A_{00,\mathrm{I}}}$ and $P_{A_{00,\mathrm{II}}}$.\\

\textit{Step 2: Fluid scaled process.} Dividing both sides of the equations by $N$ results in a fluid scaled process.
Further, because of the martingale results in \cite{bremaud1981} and \cite{pang2007martingale} we can define noise terms $e_{ij}(N)$ that tend to 0 as $N\rightarrow\infty$ with $i\ge 0$ and $j\in \{0,1\}$. The fluid scaled system can be rewritten as follows, with $i\ge2$,
\begin{equation} \label{eq:fluid_scaled}
\begin{array}{rcl}
\frac{\overline{Q}^N_{00}(t)}{N} & = & \frac{\overline{Q}^N_{00}(0)}{N} + \m_2 \intt \left(\frac{\overline{Q}_{01}^N(s)}{N} - \frac{\overline{Q}_{11}^N(s)}{N} \right) \,ds - \lam  \intt \left[ p_{00}^N(\mathrm{I},s) +  p_{00}^N(\mathrm{II},s)\right]\,ds + e_{00}(N)\\

\frac{\overline{Q}^N_{01}(t)}{N} & = & \frac{\overline{Q}^N_{01}(0)}{N} -  \m_2 \intt  \left(\frac{\overline{Q}_{01}^N(s)}{N} - \frac{\overline{Q}_{11}^N(s)}{N} \right) \,ds +\lam  \intt  p_{00}^N(\mathrm{II},s)\,ds  + e_{01}(N)\\

\frac{\overline{Q}^N_{10}(t)}{N} & = & \frac{\overline{Q}^N_{10}(0)}{N} -  \m_1 \intt  \left(\frac{\overline{Q}_{10}^N(s)}{N} - \frac{\overline{Q}_{20}^N(s)}{N} \right)\,ds + \lam  \intt  p_{00}^N(\mathrm{I},s)\,ds  + e_{10}(N)\\

\frac{\overline{Q}^N_{11}(t)}{N} & = & \frac{\overline{Q}^N_{11}(0)}{N} -  \m_1 \intt  \left(\frac{\overline{Q}_{11}^N(s)}{N} - \frac{\overline{Q}_{21}^N(s)}{N} \right)\,ds + \lam  \intt  p_{01}^N(\mathrm{I},s)\,ds  + e_{11}(N)\\

\frac{\overline{Q}^N_{i0}(t)}{N} & = &\frac{\overline{Q}^N_{i0}(0)}{N} -  \m_1 \intt  \left(\frac{\overline{Q}_{i0}^N(s)}{N} - \frac{\overline{Q}_{i+1,0}^N(s)}{N} \right)\,ds + \lam \intt  p_{i-1,0}^N(\mathrm{I},s)\,ds  + e_{i0}(N)\\

\frac{\overline{Q}^N_{i1}(t)}{N} & = & \frac{\overline{Q}^N_{i1}(0)}{N} -  \m_1 \intt  \left(\frac{\overline{Q}_{i1}^N(s)}{N} - \frac{\overline{Q}_{i+1,1}^N(s)}{N} \right)\,ds + \lam  \intt  p_{i-1,1}^N(\mathrm{I},s)\,ds + e_{i1}(N). \\

\end{array}
\end{equation}

\textit{Step 3: Towards fluid limits.} While making the transition from integral equations to differential equations with $N$ tending to infinity, the representation of the departure terms in~(\ref{eq:fluid_limit}) is straightforward. The arrival terms in the differential equations, on the other hand, are not immediately obvious.\\
To illustrate the difficulty, assume there are among the $N$ servers only a small number of idle servers. As the allocation strategy describes, one of these servers will be selected by an arriving job. If the number of idle servers is small and the arrival rate is sufficiently high, rapid switches will occur in the indicator function $\mathds{1}\{Q^N_{00}  = 0 \}$. A server that becomes idle due to a service completion will immediately be selected again by the arriving job. However, the fraction of empty servers ($Q^N_{00}/N$) will be more robust against these changes due to the fluid scaling.\\
In general, this phenomenon is called `separation of time scales' as described by Hunt and Kurtz~\cite{hunt1994large}. One observes the interaction of two processes. One process evolves very fast, namely the number of empty servers, while the second process evolves much slower, the occupancy fractions in this setting. In order to obtain the arrival terms of the fluid limit, we should be able to combine these processes. Focusing on the first arrival integral in~(\ref{eq:fluid_scaled}), the question arises how to handle the expression
\begin{equation}
\lim_{N\rightarrow\infty} \lam \intt \left[ p_{00}^N(\mathrm{I},s) +  p_{00}^N(\mathrm{II},s)\right]\,ds = \lim_{N\rightarrow\infty} \lam \intt  \mathds{1}\{Q^N_{00}(s)  > 0 \}\,ds?
\end{equation}
A similar problem is analyzed in \cite{hunt1994large} where one needs to take the limit of a integral of an indicator function. The existence of a measure $\alpha$ is deduced such that 
\begin{equation}
\lim_{N\rightarrow\infty} \lam \intt  \mathds{1}\{Q^N_{00}(s)  > 0 \}\,ds =  \lam \intt\alpha(s)\,ds.
\end{equation}
The existence of this function $\alpha$, which does not need to be continuous, can be justified by the following reasoning. In a small time interval, say $[0,\dt]$, the number of idle servers is a heavily fluctuating process, though the process describing the occupancy proportions is approximately constant. During this small interval, the number of idle servers can be considered as a birth-and-death process with `death' rate $\lam $, since an arriving job causes a reduction in the number of idle servers. The `birth' rate is determined by the occupancy proportions, i.e. the proportion of servers that are working on type~I or type~II jobs. Then it is argued in \cite{hunt1994large} that 
\begin{equation}
\frac{1}{\dt}\int_{0}^{\dt} \mathds{1}\{Q^N_{00}(s)  > 0 \}\,ds,
\end{equation}
after application of the ergodic theorem, converges to an invariant measure if $N$ tends to infinity. This invariant measure will give rise to the function $\alpha$. One already senses that the presence or absence of idle servers should be handled as two different cases. Therefore we make a distinction between $q_{00}$ strictly positive or equal to zero in the intuitive explanation of the structure of the fluid limit.\\
\textit{The case $q_{00}>0$.} When the number of idle servers is sufficiently large, each arriving job will be allocated to an idle server for sure. A fraction
\begin{equation}
\left(1-q_{00} \right)^{d_1}
\end{equation}
of the arriving jobs will be allocated as type~II jobs  which causes the changes in~(\ref{eq:fluid_limit}) for $\overline{q}_{00}$, $\overline{q}_{01}$ and $\overline{q}_{10}$.\\ 
\textit{The case $q_{00}=0$.} Idle servers are generated at rate $\m_1 q_{10} + \m_2 q_{01}$. Since $d_1$ is finite, the probability that the $d_1$-selection would contain an idle server is negligible, each idle server will be provided with a type~II job when the arrival rate is high enough. If $\tilde{\lam} = (\lam- \m_1 q_{10} + \m_2 q_{01})^+$ is strictly larger than zero, a fraction 
\begin{equation}
\frac{\m_1q_{10}+ \m_2q_{01}}{\lam} = \frac{\lam - \tilde{\lam}}{\lam}
\end{equation} 
of the stream of incoming jobs will immediately be redirected to the idle servers as a type~II job. The excess stream of incoming jobs (fraction $\tilde{\lam}/\lam$) will not observe any idle server and will start to form (type~I) queues in front of the servers of the $d_1$-selection according to a straightforward generalization of the transition probabilities mentioned in step 1.

This concludes the derivation of the fluid limit~(\ref{eq:fluid_limit}).

\subsubsection{Proof of Theorem~\ref{th:fixedpoint}: fixed points}\label{subsec:fluid_fp}
We will start with the proof of the closed-form fixed point and show that this is the only fixed point without idle servers on fluid level, i.e. $q_{00}=0$. Next, we will consider fixed points with $q_{00}>0$.\\

\textit{Fixed points with $q_{00} = 0$.} The correctness of the expression in~(\ref{eq:fixedpoint}) can easily be confirmed by substitution into~(\ref{eq:fluid_limit}). %
The result can be established in two steps. First, we observe that the derivatives of $(\overline{q}_{i0})_i$ in~(\ref{eq:fluid_limit}) remain zero once $(\overline{q}_{i0}^*)_i$ equals zero. Then, we substitute $(\overline{q}_{i0}^*)_i = 0$ into the derivatives of $(\overline{q}_{i1})_i$. For $i\ge1$ we obtain: \\
\begin{equation}\label{eq:vereenvoudiging}
\frac{d}{dt} \overline{q}_{i1}^* = \m_1(\overline{q}_{i+1,1}^* - \overline{q}_{i1}^*) + \tilde{\lam } \left[ (\overline{q}_{i-1,1}^*)^{d_1} -(\overline{q}_{i1}^*)^{d_1}\right] = 0.
\end{equation}
These equations can be solved and one obtains the fixed point as given in~(\ref{eq:fixedpoint}). Note the similarity between~(\ref{eq:vereenvoudiging}) and the fluid limit of a JSQ($d_1$) policy with reduced arrival rate 
\begin{equation}
\tilde{\lam} = \lam - \frac{\m_1-\lam}{\m_1-\m_2}\m_2 =  \lam - \m_2q_{01}^*,
\end{equation}
in a setting where each of the exchangeable servers works at rate $\m_1$ \cite{mitzenmacher2001power}.

Second, this fixed point is unique under the condition that $q_{00}$ equals zero.
From Lemma~2 in \cite{mitzenmacher2001power} we know that the fixed point of the fluid limit in the JSQ($d_1$) setting is unique when $d_1\ge 2$. This implies that under the condition that all servers have a type~II job, i.e. $\overline{q}_{i0}^* = 0$ for all $i$, uniqueness is guaranteed. Assume by contradiction that another fixed point exists without idle servers but with possibly a positive cumulative fraction $\overline{q}_{i0}^*$ for some $i$. We focus on the differential equations of $(\overline{q}_{i0})_{i\ge1}$ under this fixed point. From
\begin{equation}
\frac{d}{dt}\overline{q}_{10}^* = \m_1(\overline{q}_{20}^*-\overline{q}_{10}^*)=0
\end{equation}
we get that $\overline{q}_{10}^* = \overline{q}_{20}^*$. Repeating this procedure for $i = 2$,
\begin{equation}
\begin{array}{rcl}
\frac{d}{dt} \overline{q}_{20}^* &= &\m_1(\overline{q}_{30}^*-\overline{q}_{20}^*) + \tilde{\lam}  \left[ (\overline{q}_{10}^*-\overline{q}_{11}^*)^{d_1} -(\overline{q}_{20}^*-\overline{q}_{11}^*)^{d_1}\right]\\
%&= &\m_1(\overline{q}_{30}^*-\overline{q}_{20}^*) + \tilde{\lam}  \left[ (\overline{q}_{20}^*-\overline{q}_{11}^*)^{d_1} -(\overline{q}_{20}^*-\overline{q}_{11}^*)^{d_1}\right]\\
&= &\m_1(\overline{q}_{30}^*-\overline{q}_{20}^*)=0, \\

\end{array}
\end{equation}
results in $\overline{q}_{20}^* = \overline{q}_{30}^*$.
By induction we could show that $\overline{q}_{i0}^* = \overline{q}_{i+1,0}^*$ for $i\ge1$, this leads to $\overline{q}_{i0}^* = 0$ for $i\ge 1$. This proves the uniqueness of the fixed point when $q_{00}$ equals zero.\\

\textit{Fixed points with $q_{00} > 0$.} Under this setting, the fluid-limit equations~(\ref{eq:fluid_limit}) simplify significantly.

\begin{equation} 
\begin{dcases}
\tfrac{d}{dt} \overline{q}_{00}=\m_2(\overline{q}_{01} - \overline{q}_{11}) - \lam(1-q_{00})^{d_1} \\
\tfrac{d}{dt} \overline{q}_{01}=\m_2(\overline{q}_{11}-\overline{q}_{01}) + \lam(1-q_{00})^{d_1}\\
\tfrac{d}{dt} \overline{q}_{10}= \m_1(\overline{q}_{20} - \overline{q}_{10}) + \lam\left(1-(1-q_{00})^{d_1}\right)\\
\tfrac{d}{dt}  \overline{q}_{11}=\m_1(\overline{q}_{21} - \overline{q}_{11})\\
\text{for~}i\ge 2,\\
\tfrac{d}{dt}\overline{q}_{i0}=\m_1(\overline{q}_{i+1,0} - \overline{q}_{i0})\\
\tfrac{d}{dt} \overline{q}_{i1}=\m_1(\overline{q}_{i+1,1} - \overline{q}_{i1})
\end{dcases}
\end{equation}

For any fixed point it should hold that
\begin{equation}
\left\{
\begin{array}{rcl}
\overline{q}_{i0}^* & = & \overline{q}_{i+1,0}^* \\
\overline{q}_{i1}^* & = & \overline{q}_{i+1,1}^* 
\end{array}
\right.
\end{equation}
for $i\ge 2$, then it follows that $(\overline{q}_{i0}^*)_{i\ge 2} = 0$ and $(\overline{q}_{i1}^*)_{i\ge 1} = 0$. This implies that the only positive fractions are $q_{00}$, $q_{01}$ and $q_{10}$. Rewriting the fluid limit in a non-cumulative expression gives us:
\begin{equation} \label{eq:fllim_pos}
\begin{dcases}
\tfrac{d}{dt} q_{00}=\m_1q_{10} + \m_2q_{01}-\lam\\
\tfrac{d}{dt} q_{01}= -\m_2q_{01} + \lam(1-q_{00})^{d_1} \\
\tfrac{d}{dt} q_{10}= -\m_1q_{10} + \lam\left(1-(1-q_{00})^{d_1}\right).
\end{dcases}
\end{equation}
From the second and third equality it is clear that once $q_{00}^*$ is known, we know the entire fixed point:
\begin{equation}
\left\{
\begin{array}{rcl}
q_{01}^* & =&  \frac{\lam}{\m_2}(1-q_{00}^*)^{d_1}\\
q_{10}^* &= & \frac{\lam}{\m_1}  \left(1-(1-q_{00}^*)^{d_1}\right).
\end{array}
\right.
\end{equation}
The system in~(\ref{eq:fllim_pos}) is linearly dependent. We use the fact that $q_{00}$, $q_{01}$ and $q_{10}$ must sum up to one to determine $q_{00}$. It must hold:
\[
1 = q_{00} + (1-q_{00})^{d_1} \left( \frac{\lam}{\m_2} - \frac{\lam}{\m_1} \right) + \frac{\lam}{\m_1}.
\]
Define $x\doteq 1-q_{00}$. We are interested in the zero points of the polynomial $f$ within $[0,1)$ with
\begin{equation}
f(x) = x^{d_1} \left( \frac{\lam}{\m_2} - \frac{\lam}{\m_1} \right) - x + \frac{\lam}{\m_1}.
\end{equation}
We will evaluate the existence of the fixed points based on the behaviour of $f$ and its derivative,
\begin{equation}
f'(x) = d_1 \left( \frac{\lam}{\m_2} - \frac{\lam}{\m_1} \right) x^{d_1-1} -1 .
\end{equation}
Furthermore,
\begin{equation}
\begin{array}{rcl}
f(0)& =& \frac{\lam}{\m_1} >0 \\
f(1) &=& \frac{\lam}{\m_2}-1 >0,
\end{array}
\end{equation}
and $f'$ is monotone increasing on $(0,1)$ with
\begin{equation}
\begin{array}{rcl}
f'(0)& =& -1<0 \\
f'(1) &=& d_1\lam \left(\frac{1}{\m_2}-\frac{1}{\m_1}\right)-1.
\end{array}
\end{equation}
Since $f$ is positive in both its endpoints and the derivative $f'$ is monotone increasing, we need at least a vanishing derivative in $(0,1)$ in order to have a fixed point. This is guaranteed when $f'(1) > 0$, this is the first condition from~(\ref{eq:fixedpoint_restrict}). We now know that $f$ attains a local minimum at 
\begin{equation} \label{eq:xtilde}
\tilde{x}  \doteq \left(\frac{1}{d_1} \frac{1}{\lam \left( \frac{1}{\m_2}-\frac{1}{\m_1}\right)}\right)^{\frac{1}{d_1-1}}
\end{equation}
and is strictly positive in its endpoints. If $f(\tilde{x})$ is exactly zero, we have one fixed point, namely $q_{00}^* = 1-\tilde{x}$. But only in very special cases the second condition of~(\ref{eq:fixedpoint_restrict}) is satisfied with equality for a random choice of $d_1$, $\lam$, $\m_1$ and $\m_2$. On the other hand, if $f(\tilde{x})<0$, i.e. if  also
\begin{equation}
\left(1- \frac{1}{d_1}\right) \left( \frac{1}{d_1}\right)^{\frac{1}{d_1-1}}  >  \frac{\lam}{\m_1} \left( \lam \left( \frac{1}{\m_2} - \frac{1}{\m_1}\right) \right) ^{\frac{1}{d_1-1}}
\end{equation}
holds, we have exactly two fixed points such that $q_{00} + q_{01}+ q_{10} = 1$. There is one fixed point situated at each side of $\tilde{x}$ in the interval $(0,1)$. This gives that for $d_1$ large enough we can find two solutions of the reduced system of differential equations. It can be shown by contradiction that both fixed points are larger than $\lam/\m_1$, so the corresponding fractions of idle servers is smaller than $1-\lam/\m_1$.\\

For completeness we mention that $\lam = \m_2$ would imply that $f(1) = 0$ and so the proportion of empty servers is zero which violates the assumption that $q_{00}>0$. Moreover, if $\lam < \m_2$, then the polynomial $f$ vanishes in the interval $(0,1)$. The monotone increasing property of the derivative of $f$ leads to the fact that there exists a unique fixed point $x^*$ in $(0,1)$. This results in a unique fixed point $(q_{00}^*, q_{01}^*,q_{10}^*)$ with $q_{00}^*>0$.

This concludes the proof of Theorem~\ref{th:fixedpoint}.

%%%%%%%%%%%%%%%%%%%%%%%%%%%%%%%%%%%%%%%%%%%%%%%%%%%%%%%%%%%%%%%%%%%%%%%%%%%%%
\subsubsection{Proof of Theorem~\ref{th:fixed_unst}: local \textup{(}in\textup{)}stability} \label{subsec:proof_loc_stab}

We will prove local (in)stability using the indirect Lyapunov method based on the Hartman-Grobman Theorem \cite{hartman1960local}. This theorem states that a system of differential equations behaves near its fixed points as its linearized version. The eigenvalues of the linearized system will define the local behavior of the system unless one of the eigenvalues has a real part equal to zero, then the Hartman-Grobman theorem is inconclusive. When we would immediately apply this theorem to one of the two fixed points of~(\ref{eq:fllim_pos}) we obtain an eigenvalue exactly equal to zero, but one can resolve this issue since~(\ref{eq:fllim_pos}) is a redundant system. Since $q_{00} + q_{01} + q_{10} = 1$, it is sufficient to know the instantaneous change of two variables. Each elimination will lead to the same two eigenvalues so we can remove for instance the third equation from~(\ref{eq:fllim_pos}):
\begin{equation}
%\begin{array}{rcl}
\begin{dcases}
\tfrac{d}{dt}q_{00}  =  \m_1(1-q_{00}-q_{01})+ \m_2q_{01} - \lam \\
\tfrac{d}{dt}q_{01}  =  -\m_2q_{01}  + \lam(1-q_{00})^{d_1}.
\end{dcases}
%\end{array}
\end{equation}

Let $(q_{00}^*,q_{01}^*,q_{10}^*)$ denote a fixed point, then the matrix of the linearized system looks as follows near its fixed point:
\begin{equation}
  \begin{bmatrix}
-\m_1 & \m_2-\m_1 \\ -\lam d_1(1-q_{00}^*)^{d_1-1} & -\m_2
\end{bmatrix}.
\end{equation}
The corresponding eigenvalues are given by
\begin{equation}
\alpha_{\pm} = \frac{1}{2}\left[ -(\m_1+\m_2) \pm \sqrt{(\m_1-\m_2)^2+4\lam d_1(\m_1-\m_2)(1-q_{00}^*)^{d_1-1}}\right].
\end{equation}
Since $\m_1 > \mu_2$ the quantity under the root is always positive, so the square root is real. This implies furthermore that $\alpha_- <0$. To determine the sign of $\alpha_+$ we need to make a distinction between the two fixed points. From the proof of Theorem~\ref{th:fixedpoint}, we know that the two fixed points are on both sides of $\tilde{x}$, with $\tilde{x}$ as in~(\ref{eq:xtilde}). For
\begin{equation}
1-q_{00}^* > \tilde{x} = \left( \frac{1}{d_1}\frac{\m_1\m_2}{\lam (\m_1-\m_2)}\right)^{1/(d_1-1)}
\end{equation}
we have that
\begin{equation}
\begin{array}{rcl}
2\alpha_+ & > &  -(\m_1+\m_2) + \sqrt{(\m_1-\m_2)^2+4\lam d_1(\m_1-\m_2)( \frac{1}{d_1}\frac{\m_1\m_2}{\lam (\m_1-\m_2)})}\\
& =&  -(\m_1+\m_2) + \sqrt{(\m_1+\m_2)^2} \\
& =& 0.
\end{array}
\end{equation}
This shows that the fixed point with the smallest  proportion of idle servers is unstable.

 When $1-q_{00}^* < \tilde{x}$, it follows in a similar way that
$
2\alpha_+  < 0.
$
This shows that the fixed point with the largest proportion of idle servers is locally stable.
This concludes the proof of Theorem~\ref{th:fixed_unst}.

%%%%%%%%%%%%%%%%%%%%%%%%%%%%%%%%%%%%%%%%%%%%%%%%%%%%%%%%%%%%%%%%%%%%%%%%%%%%%%%%%%%%%%%%%%%%%%%%%%%%%%%%%%%%%%%%%%%%%%%%%%%%%%%%%%%%%%%%%%%%%%%%%%%%%%%%%%
\section{Conclusion and outlook} \label{sec:conclusion}
We investigated load balancing issues in a service system where particular servers are better equipped to process certain jobs due to affinity or compatibility relations.
The general model in particular covers the setting with an underlying network topology~$G_N$, referred to as the \textit{graph model}. The analysis of the graph model is severely complicated by the lack of exchangeability among the servers, a feature linked to the supermarket modeling framework that allows mean-field techniques. We constructed the a novel \textit{affinity coupling} to obtain a stochastic performance bounds for the general model and more specific settings, for instance model instances where the underlying graph topology~$G_N$ has a specific minimum degree or is a $d$-regular graph.

Another instance of the general model, the \textit{combinatorial model}, has enough inherent symmetry to conduct a fluid-limit analysis. The fluid limit was stated in terms of a set of discontinuous differential equations and its fixed point sensitively depends on the size~$d$ of the primary selection.
 When $d$ is sufficiently small, a unique fixed point exist but the associated waiting time does not vanish. 
When the primary selection is sufficiently large, a fixed point arises that does provide a zero waiting time. On the other hand, the above-mentioned fixed point still persists, giving rise to bistability issues. 
 
 As mentioned above, the stochastic upper bounds for the graph model in terms of a supermarket model with a JSQ($d$) policy require the degrees in the underlying graph to be relatively high compared to~$d$.
To some extent, this reflects that the performance may be poor in certain pathological cases even when the node degrees are fairly high. An interesting topic for further research would be to extend the affinity coupling and possibly identify relevant structural conditions on the graph topology in order to sharpen these bounds.

Recall that a supermarket model with a JSQ($d$) policy is equivalent to the combinatorial model with server selections of size~$d$ and identical arrival rates when jobs could not be allocated as a type~II job.
A natural conjecture is that the latter combinatorial model is the best-case scenario given a maximum cardinality~$d$ of the server selections.
This would imply that the supermarket model with a JSQ($d$) policy provides stochastic lower bounds in some appropriate sense for \textit{any} affinity model with server selections of size at most~$d$.

The bistability of the fluid limit of the combinatorial model for large values of~$d_1$ not only precludes any convergence statements for the stationary distribution, but also suggests that the allocation strategy could possibly be refined.
In future work we intend to examine such refinements and establish that these eliminate the \textit{queueing} fixed point and render the \textit{no-queueing} fixed point globally stable.

%\begin{acknowledgements}
%\end{acknowledgements}

\bibliographystyle{abbrv}
\bibliography{references}   % name your BibTeX data base

\end{document}